%% file: thesis_teste.tex
\documentclass[
a4paper,
12pt,
times, % DON'T USE THIS OPTION
numbered,
print,
index,
oneside, % just for shortening the example
]{PhDThesisPSnPDF}

\usepackage{hyperref}
\usepackage{enumitem}
\usepackage[utf8]{inputenc}
\usepackage{stmaryrd}
\usepackage{dsfont}
\usepackage{amsmath}
\usepackage{amsfonts}
\usepackage{mathtools}

\usepackage[T1]{fontenc} 
\usepackage{amssymb}
\usepackage{relsize,exscale}
\input{Preamble/preamble}
\newtheorem{theoreme}{ Th\'{e}or\`{e}me}[section]
\newtheorem{proposition}{ Proposition}[section]
\newtheorem{lemma}[theoreme]{Lemme}

\newenvironment{preuve}[1][Preuve]{\begin{trivlist}
		\item[\hskip \labelsep {\bfseries #1}]}{\end{trivlist}}
\newenvironment{defini}[1][Définition]{\begin{trivlist}
		\item[\hskip \labelsep {\bfseries #1}]}{\end{trivlist}}

\newenvironment{remark}[1][Remarque]{\begin{trivlist}
		\item[\hskip \labelsep {\bfseries #1}]}{\end{trivlist}}

\newcommand{\qed}{\nobreak \ifvmode \relax \else
	\ifdim\lastskip<1.5em \hskip-\lastskip
	\hskip1.5em plus0em minus0.5em \fi \nobreak
	\vrule height0.75em width0.5em depth0.25em\fi}

\input{thesis-info}

\ifdefineAbstract
\fi

\begin{document}

	\frontmatter
	
	\maketitle
	
	%\include{Dedication/dedication}
	\include{Declaration/declaration}
	\include{Acknowledgement/acknowledgement}
	\include{Abstract/abstract}

	% *********************** Adding TOC and List of Figures ***********************
	
	\tableofcontents
	
	%\listoffigures
	
	%\listoftables
	
	% \printnomenclature[space] space can be set as 2em between symbol and description
	%\printnomenclature[3em]
	
	\printnomenclature
	
	% ******************************** Main Matter *********************************
	\mainmatter
	
	\include{Chapter1/chapter1}
	\include{Chapter2/chapter2}

	\include{Chapter3/chapter3}

\begin{spacing}{0.9}
	
	% To use the conventional natbib style referencing
	% Bibliography style previews: http://nodonn.tipido.net/bibstyle.php
	% Reference styles: http://sites.stat.psu.edu/~surajit/present/bib.htm
	
	\bibliographystyle{apalike}
	\bibliographystyle{plainnat} % use this to have URLs listed in References
	\cleardoublepage
	\bibliography{References/references.bib} % Path to your References.bib file

	% If you would like to use BibLaTeX for your references, pass `custombib' as
	% an option in the document class. The location of 'reference.bib' should be
	% specified in the preamble.tex file in the custombib section.
	% Comment out the lines related to natbib above and uncomment the following line.
	
	%%%%%%%%%%%%%%%\printbibliography[heading=bibintoc, title={References}]

\end{spacing}

\end{document}

%% file: thesis-info.tex
\title{Autour de la fonction Zêta de Riemann}
\subtitle{Maître de stage : Yang Hui He \newline Tuteur école \, \, : Thibaut LEGOUIC } 
\author{EL HASSANI Ansar}

\dept{Department of Mathematics}

\university{City University of London}

\degreetitle{Stage S8 recherche international}

\college{Ecole Centrale Marseille}

%% file: Declaration/declaration.tex
% ******************************* Thesis Declaration ***************************

%% file: Acknowledgement/acknowledgement.tex
% ************************** Thesis Acknowledgements **************************

\begin{acknowledgements}      

J'aimerais remercier mon Maître de stage Yang Hui He pour tout le temps et les discussions qu'il m'a accordé plusieurs fois par semaines durant toutes la durée du stage . 

Je remercie mes collègues doctorants de m'avoir présenter leurs travailles , pour m'avoir expliqué quelques concepts et pour le temps qu'on a partagé ensemble .

Je remercie mes colocataires , avec qui ce séjour a Londres était une excellente expérience .

Je remercie Oumaima SADIKI , Ibtihaj SADIKI , Olivier ROUX , Mohammed ALLALI , Nizar Bellazrak et plus particulièrement Mon Père pour leurs nombreuses corrections et leurs relectures .

Je remercie Olivier Ramaré pour sa relecture ainsi que ses très nombreuses corrections et commentaires  . 

Je remercie mes parents , ma famille et mes amis pour tout .

\end{acknowledgements}

%% file: Abstract/abstract.tex
% ************************** Thesis Abstract *****************************
% Use `abstract' as an option in the document class to print only the titlepage and the abstract.
\begin{abstract}
	In the first chapter, we will present a computation of the square value of the module of L functions associated to a Dirichlet character.
	
	This computation suggests to ask if a certain ring of arithmetic multiplicative functions exists and if it is unique. This search has led to the construction of that ring in chapter two.
	
	Finally, in the third chapter, we will present some propositions associated with this ring.
	The result below is one of the main results of this work.	
	
	\[~~~~~~\]
	Dans le premier chapitre, je présente un calcul du module au carré des fonctions L associé à un caractère de Dirichlet.
	
	Ce calcul m'a poussé à explorer l'existence et l'unicité d'un certain anneau de fonctions arithmétiques multiplicatives que je construis par la suite dans le chapitre deux. 
	
	Enfin, dans le chapitre trois, je présente quelques propositions associées à cet anneau. 
	
	Un des résultats principaux est présenté ci-dessous: 
	
	\[~~~~~~\]
	
	\textnormal{For F and G two completely multiplicative functions, $ s $ a complex number such as the dirichlet series $ D(F,s) $ and $ D(G,s) $ converge :}
	\[\forall F,G \in \mathbb{M}_{c} :  D(F,s) \times D(G,s) = D(F \times G,2s) \times D( F {\, \scalebox{0.5}{$\square$} \,} G ,s) \]
	Here are some similar versions, with $ s = x+iy $ :
	\[  \forall F , G \in \mathbb{M}_{c} : ~ D(F,s) \times D(G,\overline{s}) = D(F \times G,2x) \times D( \frac{F}{\text{Id}_{e}^{iy}} \scalebox{0.6}{$\square$} \frac{G}{\text{Id}_{e}^{-iy}} , x ) \]
	\[  \forall F , G \in \mathbb{M}_{c} : ~ |D(F,s)|^{2} = D(|F|^{2},2x) \times D( \frac{F}{\text{Id}_{e}^{iy}} \scalebox{0.6}{$\square$}  \overline{\frac{F}{\text{Id}_{e}^{iy}}} , x ) \]
	
\end{abstract}

%% file: Chapter1/chapter1.tex
\chapter{Factorisation du Module au carré des fonctions L }

\label{ch:Factorisation du Module au carré des fonctions L}

\section{Présentation du problème  }
% body of thesis comes here
\subsection{La fonction zêta  }
La série zêta a été découverte par Leonhard Euler et est définie comme suit :
\[ \forall s \in \mathbb{C} , ~   \Re(s)>1 \: ~ \zeta(s)  = \sum_{n=1}^{+\infty} \frac{1}{n^{s}} \]
La série $ \sum_{n=1}^{+\infty} \frac{1}{n^{s}} $ n'est pas convergente pour $ \Re(s) \le  1 $ , néanmoins on peut définir un prolongement (unique au sens du prolongement analytique ) sur tout le plan complexe, sauf en $ s = 1 $.
\\

Ce prolongement analytique est ce qu'on appelle la fonction zêta \cite{Edwards} \cite{Riemann} .

\subsection{Le produit Eulérien  }

Euler a aussi a établi le produit "Eulérien" associé à la fonction zêta :

\[ \forall s \in \mathbb{C} , ~   \Re(s)>1  ~  \zeta(s) = \sum_{n=1}^{+\infty} \frac{1}{n^{s}} = \prod_{p \in \mathbb{P}} \frac{1}{1-\frac{1}{p^{s}}} \]
Pour comprendre plus facilement l'importance de cette formule , il est plus adéquat de l'écrire comme ceci :

\[ \forall s \in \mathbb{C} ,  ~  \Re(s)>1  ~  \prod_{p \in \mathbb{P}} \frac{1}{1-\frac{1}{p^{s}}} = \prod_{p \in \mathbb{P}} 1+\frac{1}{p^{s}}+\frac{1}{p^{2s}}+\frac{1}{p^{3s}}+... \]

\[ \forall s \in \mathbb{C} ,  ~  \Re(s)>1  ~ \zeta(s)  = \sum_{n=1}^{+\infty} \frac{1}{n^{s}} = 1+\frac{1}{2^{s}}+\frac{1}{3^{s}}+\frac{1}{4^{s}}+...  \]

L'idée est que la formule $ \sum_{n=1}^{+\infty} \frac{1}{n^{s}} $ présente les entiers en tant que suite de nombres espacés par un pas de 1 \textit{(vision additive  / métrique)} , alors que la formule $ \prod_{p \in \mathbb{P}} \frac{1}{1-\frac{1}{p^{s}}} $ exprime les entiers en tant que combinaison de différentes puissances de nombres premiers \textit{(vision multiplicative / combinatoire )}.

\subsection{Equation fonctionnelle  }

Riemann a prouvé \cite{Riemann} que la fonction zêta vérifie l'équation fonctionnelle suivante :
\[ \forall s \in \mathbb{C}-1 ~ : ~ \zeta(s) = 2^s \pi^{s-1} \sin \left( \frac{\pi s}2\right) \Gamma(1-s)\zeta(1-s) \]
Ce résultat exprime une relation entre les valeurs $ \zeta (s) $ et $ \zeta(1-s) $ , et plus particulièrement , il dit que l'ensemble des zéros de la fonction zêta est symétrique par rapport à la droite $ \Re(s) = \frac{1}{2} $ .

\subsection{L'Hypothèse de Riemann }

Bernhard Riemann est un des meilleurs mathématiciens de son époque . Il connaissait bien le domaine des transformées de Fourier et il avait réussi à lier la théorie des nombres à celle des fonctions à variables complexes . Dans un certain sens , il a trouvé que la transformée de Fourier de l'ensemble des nombres premiers , et cela dans le seul article \cite{Riemann} (très court) sur la théorie des nombres qu'il a écrit.

Il a compris que la fonction zêta avait des liens importants avec les nombres premiers , son travail a débouché sur une formule exacte qui permet de trouver l'emplacement de chaque nombre premier : 

La fonction $ li $, dîte logarithme intégral, est définie par :
\[\operatorname{li}(x) = \int_{0}^{\infty} \frac{dt}{\ln(t)}\]
$ \rho $ sont les zéros non triviaux de la fonction zéta :
\[ \forall \rho ~:~  \zeta(\rho) = 0 ~ , ~ 0 < \Re(\rho) < 1 \]
$ \psi $ est une fonction qui permet de retrouver la fonction de comptage $ \pi $ , nous n'aurons pas besoin de sa définition précise .

\[ \psi(x) = \operatorname{li}(x) - \sum_\rho \operatorname{li}(x^\rho) -\log(2) +\int_x^\infty\frac{dt}{t(t^2-1)\log(t)} \]
\\
Cette formule a malheureusement un problème important : elle contient une somme sur les zéros non triviaux $ \rho $  ( c'est à dire de partie réel $ 0 < \Re(\rho) < 1 $  )  de la fonction zêta .
L'étude des zéros de la fonction zêta prend alors beaucoup d'importance , l'hypothèse qu'avait fait Riemann était que tout les zéros avec leur partie réel égale à $ \frac{1}{2} $ . 
\subsection{Le théorème des nombres premiers }
\subsubsection{ L'ordre dans les nombres premiers }
Le théorème des nombres premiers dit : 
\[ \pi(x)\sim \operatorname{li}(x) \]
Où $ \pi $ est la fonction de comptage des nombres premiers .

l'Hypothèse de Riemann est équivalente au résultat suivant : 
\[ \pi(x) = Li(x) + O(\sqrt{x} \cdot \ln(x)) \]
Où $ \pi $ est la fonction de comptage des nombres premiers .
\\

l'Hypothèse de Riemann complète alors le Théorème des Nombres Premiers en précisant l'ampleur des fluctuations entre $ \pi(x) $ et $ li(x) $ . \\

il est aussi très intéressant de savoir que l'un des points clefs des premières preuves du théorème des nombres premiers est l'absence de zéro sur la droite $ \Re(s)=1 $ . 

Il faut se rappeler que les zéros de la fonction zêta sont liés par la transformée de Fourrier au nombres premiers , ils peuvent alors être interprétés comme les notes de la musique des nombres premiers (ceci reste une interprétation simplifier) . Si elle s'avère correcte , L'hypothèse de Riemann correspond alors à une régularité ordonné sur l'ensemble des nombres premiers .
\subsubsection{Le chaos dans les nombres premiers }

Voici un autre équivalent de l’hypothèse de Riemann , Trouvé par Landau dans sa thèse ( longue de 13 page ) \cite{Landau} :
\[ \forall \epsilon > 0 ~ : ~ \lim_{n \rightarrow \infty} \frac{\lambda(1)+\lambda(2)+\lambda(3)+...+\lambda(n)}{n^{\frac{1}{2}+\epsilon}} =0 \]

Où $ \lambda(l) = (-1)^{\omega(l)} $ et $ \omega(l) $ le nombre de nombres premiers distincts dans la décomposition de $l$   .
\\

Une première interprétation de cette limite est que les $ (-1)^{\omega(l)}  $ donnent une suite de pas dans $ \mathbb{Z} $ , le $ \frac{1}{2} $ qui apparait dans la prédiction  $ \Re(s) =  \frac{1}{2} $ de l'hypothèse de Riemann correspond alors à la racine carré obtenue en utilisant dans le théorème central limite lors du calcul de la distance à l'origine d’une marche aléatoire unidimensionnelle et isotrope pour un grand nombre de pas , un article de França et Leclair observe les mêmes phénomènes sur d'autres grandeurs liées à la  théorie des nombres \cite{LeClair8} .
\\

Cela signifierait que certaines caractéristiques ( ici $ \omega(l) $ ) des nombres premiers sont distribuées aléatoirement . Ce qui a priori est contradictoire avec la première interprétation . 
\subsubsection{les nombres premiers : l'ordre dans le chaos }
L'idée est que les deux interprétations sont complémentaires : \\

L'Hypothèse de Riemann décrit un schéma , une régularité que les nombres premiers sont censés suivre . Ce schéma c'est le fait qu'ils sont , au moins en partie , distribués aléatoirement . Le fait de savoir cette information nous permet de mieux les encadrer .
\\ 

L'article de Riemann a impressionné les mathématiciens de l'époque , particulièrement grâce aux liens qu'il met en évidence entre deux disciplines mathématiques distantes . La recherche a alors commencé pour prouver "l'Hypothèse" que Riemann a laissée , leurs efforts n'ont malheureusement pas abouti à ce jour .

\subsection{Approches de résolution à travers l'histoire  }

\cite{Connes} Une des stratégies que les mathématiciens ont appliquée pour tenter de répondre à cette question est la généralisation : 

Ils ont essayé de définir plusieurs classes de fonctions qui présentent les mêmes propriétés fondamentales de la fonction zêta (Produit Eulérien , Équation Fonctionnelle , Hypothèse de Riemann ... ) et cherché à comprendre ces nouvelles fonctions . Cela permet de ne se concentrer que sur les propriétés importantes reliées à l'Hypothèse de Riemann . 
\\

Plusieurs classes de fonctions ont ainsi été définies : 
\\

les fonctions L  , la classe de Selberg , les fonctions zêta de plusieurs variables , les fonctions zêta
de Hurwitz , les fonctions zêta sur des corps finis ...
\\

Certaines de ces classes sont au cœur de grand programme de recherche :
\\

Le programme de Langlands : ce programme essaie de relier plein d'objets mathématiques , particulièrement les formes modulaires et les fonctions L , il est important de noter qu'un des points clefs de la preuve du dernier théorème de Fermat est une correspondance qui s'inscrit dans ce programme.
\\

Un autre programme célèbre est le le programme de Polya-Hilbert dont le but est d'explorer les liens avec l'analyse fonctionnelle , avant de trouvé une interprétation en physique quantique . Un équivalent de l'hypothèse de Riemann physique (la construction d'un certain atome théorique ) a été trouvé , et plus récemment , des calculs numériques confirment une conjoncture (Conjecture de Montgomery) qui permet de prévoir plusieurs informations sur les zéros de la fonction zêta. 

\subsection{Preuve de Weil }

\cite{Hindry} Sur toutes ces classes , seule l’hypothèse de Riemann sur les corps finis a été prouvée par André Weil , cette preuve étant incomplète , elle a été par la suite complétée par toute une série de travaux de Weil mais aussi Alexandre  Grothendieck a travers son ouvrage célèbre \textit{éléments de géométrie algébrique}  et Pierre Deligne  ( en construisant tout un champ de la géométrie algébrique : schémas ,  cohomologie étale ... ) 

Ce qu'il faudrait maintenant , c'est d'arriver à reprendre les idées de la preuve de Weil , et de les appliquer au problème de base , mais il se trouve que ce n'est pas très évident .
\\

\cite{ConnesTrace} Le travaille d'Alain Connes s'inscrit dans cette optique . Durant les dernières décennies , Alain Connes a développé une théorie nommée Géométrie Non Commutative après l'étude de certains opérateurs . Il se trouve que dans cette théorie , il a réussi à donner une nouvelle preuve du fait qu'il y est une infinité de zéros sur la droite $ \Re(s)=\frac{1}{2} $ , ce résultat était prouvé par Hardy différemment. il continue de travailler jusqu'à aujourd'hui sur la question , il a notamment avancé très récemment dans la construction d'un espace qu'il espère convenable pour la compréhension des nombres premiers .

\subsection{Contexte de ce travail }

Mon tuteur avait trouvé un équivalent de la Conjoncture abc en  physique théorique  ( l’Hypothèse de Riemann , la Conjoncture abc et la Conjecture de Birch et Swinnerton-Dyer sont considérées souvent parmi les problèmes les plus importants de la théorie des nombres ) , cette trouvaille a permis de révéler un lien inattendu entre théorie des nombres et physique théorique \cite{YangHuiHe} . Ceci l'a encouragé alors à mieux comprendre l'Hypothèse de Riemann ainsi que les relations avec son domaine .
\\ 

Mon travail pendant le stage était alors simplement d'explorer l'Hypothèse de Riemann et de permettre à mon maître de stage de mieux la comprendre.

\section{Définition : }

\begin{defini}
	On définit l'ensemble des fonctions multiplicatives comme suit :
	\[
	F \in \mathbb{M} \iff 
	\begin{array}{l rcl}
	{F} \textnormal{ : }  \mathbb{N}^{\star}  \longrightarrow  \mathbb{C} \\
	F(1)=1 \\
	a \wedge b =1 \Rightarrow F(ab)=F(a)F(b) \\
	\end{array}
	\]
\end{defini}

\begin{remark}
	 $ \forall a , b \in \mathbb{N^{\star}} ~ : ~ a \wedge b $: désigne le plus grand diviseur commun de a et b 
\end{remark}

\begin{defini}
	Un caractère de Dirichlet est une fonction multiplicative qui vérifie :
	\[
	{\chi} \textnormal{ : }  \mathbb{N}  \longrightarrow  \mathbb{C} 
	\]
	\[\forall a,b \in \mathbb{N} \Rightarrow \chi(ab)=\chi(a)\chi(b)\]
	\[\exists k \in \mathbb{N^{\star}} , \forall n \in \mathbb{N} ~ \chi(n) = \chi(n + k) \]
	\[ \forall n \in \mathbb{N^{\star}} ~ n \wedge k > 1 \Rightarrow \chi(n) = 0
	\]
\end{defini}

\begin{remark}
	Un caractère de Dirichlet définit une fonction multiplicative.
\end{remark}

\begin{defini}
	La série de Dirichlet génératrice associée à une fonction multiplicative est :
	\[\forall F \in \mathbb{M} , \forall s \in \mathbb{C} , \Re(s) > \sigma  : ~ D(F,s)=\sum_{n=1}^{\infty} \frac{F(n)}{n^{s}} \]
	$\sigma $ définit une région de $\mathbb{C}$ où la série converge.

\end{defini}

\begin{defini}
	Les séries L sont les séries de Dirichlet génératrice associées à un caractère de Dirichlet :
	\[\forall s \in \mathbb{C} , \Re(s)>1 : ~ L(s,\chi) = \sum_{n=1}^{\infty} \frac{\chi(n)}{n^{s}} \]
	
\end{defini}

\begin{defini}
	Les fonctions L sont le prolongement analytique des séries L :
	
\end{defini}

\begin{defini}
	L'hypothèse de Riemann généralisée (GRH) est l'assertion que la proposition suivante est vraie pour toutes les fonctions L :
	
	\[ \forall s \in \mathbb{C} ~,~ 0<\Re(s)<1 : L(\chi,s)=0 \Rightarrow \Re(s)=\frac{1}{2} \]
	
	La fonction zêta est un cas particulier des fonctions L .
\end{defini}

\begin{remark}
	L'ensemble des nombres complexes tel que $  0<\Re(s)<1   $ est souvent appelé la bande critique .
	
	La droite des nombres complexes définie par $ \Re(s)=\frac{1}{2} $ est appelée la ligne critique .
\end{remark}

\begin{remark}
	On va étudier les fonctions L de manière générale , le cas de zêta en découle naturellement .
\end{remark}

\section{Démarche :}

Au début de mon stage , je cherchais à trouver parmi les propriétés de la fonction $ \zeta $ celles qui faisaient apparaître la droite $ \Re(z) = \frac{1}{2} $ de manière particulière , un peu comme l'équation fonctionnelle le fait en affirmant que c'est une droite de symétrie de l'ensemble des zéros de $ \zeta $. \\

Après quelque temps à travailler sur le problème , je pensais de plus que la partie imaginaire ne changeait pas beaucoup le problème .

Par exemple , si on réécrit le terme de cette série de la manière suivante : 
\[ \forall s \in \mathbb{C}-{1} , \Re(s) > 0 ~ : ~ \zeta(s) = \frac{1}{1-2^{1-s}}\cdot \sum_{n=1}^{\infty} \frac{(-1)^{n-1}}{n^{s}} = \frac{1}{1-2^{1-s}}\cdot \sum_{n=1}^{\infty} \frac{1}{n^{x}} \cdot e^{i \cdot (-yln(n))} \cdot (-1)^{n-1}   \]

On a alors :

\[|e^{i \cdot (-yln(n))} \cdot (-1)^{n-1}| = 1 \]

le $ \frac{1}{2} $ peut alors représenter une vitesse de convergence .

Après cela ,  j'ai observé que les termes indépendants de la partie imaginaire divergeaient pour $ \Re(s)=\frac{1}{2} $ dans le calcul de $ |\zeta(s)|^{2} $ avec la formule 

\[ \forall s \in \mathbb{C} ,  ~  \Re(s)>1  ~|\zeta(s)|^{2}  = \sum_{n=1}^{+\infty}\sum_{q=1}^{+\infty} \frac{1}{n^{s}} \cdot \frac{1}{q^{\bar{s}}} = \zeta(2 \cdot \Re(s)) + \sum_{n,q=1 ~ n \ne q}^{+\infty} \frac{1}{n^{s}} \cdot \frac{1}{q^{\bar{s}}} \]

Ce qui m'a guidé vers une factorisation des séries L .

\section{Somme selon les couples premiers entre eux : }

\begin{lemma}
1	\textnormal{ Soit $ \chi $ un caractère de Dirichlet, soit $ L( . ,\chi) $  sa série L  associée , soit $ {L}_{N}( . ,\chi) $ avec $ N \in \mathbb{N}^{\star} $ la somme partielle de la série L. le module de la somme partielle de la série L: }
	\[ |L_{N}(z,\chi)|^{2} = \sum_{n=1}^{N} \sum_{q=1}^{N} \frac{\chi(n)}{n^{z}}
	\overline{\frac{\chi(q)}{q^{z}}}\]
\end{lemma}

\begin{preuve}
	\begin{align*}
	|L_{N}(z,\chi)|^{2}
	& = \Bigg[ \sum_{n=1}^{N} \frac{\chi(n)}{n^{z}} \Bigg] \times \Bigg[ \sum_{q=1}^{N} \overline{\frac{\chi(q)}{q^{z}}} \Bigg] \\
	& = \sum_{n=1}^{N} \sum_{q=1}^{N} \frac{\chi(n)}{n^{z}}
	\overline{\frac{\chi(q)}{q^{z}}} 
	\end{align*}
\end{preuve}

\begin{lemma}
	\textnormal{ pour tout $ a,b \in \mathbb{N}^{*} $ , l'\'{e}criture : }
	\[ \left\{
	\begin{array}{ll}
	(a,b)=(pn, pq) \\
	n \wedge q = 1 \\
	n,q,p \in \mathbb{N}^{*}
	\end{array} \right.
	\]
	\textnormal{ existe et elle est unique }
\end{lemma}

\begin{preuve}
	soit $ n \wedge q = 1 $ et $ e \wedge f = 1 $ avec $ n,q,e,f \in \mathbb{N}^{*} $, soit  $p,l \in \mathbb{N}^{*} $ :
	
	\[ (pn, pq) = (le, lf) \Leftrightarrow \left\{
	\begin{array}{ll}
	pn = le  \\
	pq = lf
	\end{array}
	\right. \Rightarrow \frac{n}{q} = \frac{e}{f} \Rightarrow nf = qe \]
	
	en utilisant le lemme de gauss :
	
	$ nf = qe $ et $ n \wedge q = 1 $ donc $  n|e $
	
	$ nf = qe $ et $ e \wedge f = 1 $ donc $  e|n $
	
	d'ou : $ n=e $
	
	donc : $ (n, q) = (e, f) $
	
	d'ou : $ p = l $
	
	pour tout  $ a,b \in \mathbb{N}^{*} $ : 
	\[ \frac{a}{a \wedge b} \wedge \frac{b}{a \wedge b} = 1\]
	et :
	\[ \frac{a}{a \wedge b},\frac{b}{a \wedge b},{a \wedge b} \in \mathbb{N}^{*}\]
	donc l'\'{e}criture sous la forme : 
	\[ \left\{
	\begin{array}{ll}
	(a,b)=(pn, pq) \\
	n \wedge q = 1 \\
	n,q,p \in \mathbb{N}^{*}
	\end{array}
	\right.
	\]
	existe pour tout $ a,b \in \mathbb{N}^{*} $ et est unique .
\end{preuve}
\begin{lemma}
	\[|L_{N}(z,\chi)|^{2}=\sum_{n=1}^{N}\sum_{q=1}^{N} \frac{\chi(n)}{n^{z}}
	\overline{\frac{\chi(q)}{q^{z}}} = \sum_{\underset{n\wedge q=1}{n=1 ~ q=1}}^{N} 
	\sum_{{p=1}}^{\lfloor\min(\frac{N}{n},\frac{N}{q})\rfloor} \frac{\chi(np)}{(np)^{z}}
	\overline{\frac{\chi(qp)}{(qp)^{z}}} \]
\end{lemma}

\begin{preuve}
	
	la preuve consiste à combiner les deux lemmes pr\'{e}c\'{e}dents en sommant la somme du lemme 1 avec la partition du lemme 2 :
	
	\begin{align*}		
	\begin{split}
	\begin{cases}
	\begin{array}{ll}
	1 \le n \le N \\
	1 \le q \le N \\
	n,q,N \in \mathbb{N}^{*}
	\end{array}
	\end{cases}
	& \Leftrightarrow  
	\begin{cases}
	\begin{array}{ll}
	1 \le ap \le N \\
	1 \le bp \le N \\
	a \wedge b = 1 \\
	a,b,p,N \in \mathbb{N}^{*}
	\end{array}
	\end{cases}
	\\ &\Leftrightarrow  
	\begin{cases}
	\begin{array}{ll}
	\frac{1}{a} \le p \le \frac{N}{a} \\
	\frac{1}{b} \le p \le \frac{N}{b} \\
	a \wedge b = 1 \\
	a,b,p,N \in \mathbb{N}^{*}
	\end{array}
	\end{cases}
	\\ &\Leftrightarrow  
	\begin{cases}
	\begin{array}{ll}
	\max(\frac{1}{a},\frac{1}{b}) \le p \le
	\min(\frac{N}{a},\frac{N}{b}) \\
	a \wedge b = 1 \\
	a,b,p,N \in \mathbb{N}^{*}
	\end{array}
	\end{cases}
	\\ &\Leftrightarrow  
	\begin{cases}
	\begin{array}{ll}
	1 \le p \le \min(\frac{N}{a},\frac{N}{b})  \\ 
	a \wedge b = 1 \\
	a,b,p,N \in \mathbb{N}^{*}
	\end{array}     
	\end{cases}
	\end{split}
	\end{align*}
	d'où la somme demandé .
	$ \square $
\end{preuve}

\begin{remark}
	Les séries présentées ne convergent absolument que pour $ \Re(s) > 1 $ , la limite de la somme partielle ne pourra alors être prise que dans ce cas . Le reste de ce travail sera effectuer pour $ \Re(s) > 1 $ .
\end{remark}

\section{Factorisation du module au carré pour $ \Re(s) > 1 $ :}

\begin{theoreme}
	\[ |L(z,\chi)|^{2} = \Bigg[ \sum_{p=1}^{\infty} \frac{|\chi(p)|^{2}}{p^{2x}} \Bigg] \Bigg[ \sum_{l=1}^{\infty} \frac{1}{l^{x}}  \Bigg[ \delta_{l=1}(l) + \sum_{\underset{mk=l ~ m \le k }{m \wedge k = 1} } 2\Re\Big(\frac{\chi(m)}{{m}^{iy}}\overline{\frac{\chi(k)}{{k}^{iy}}}\Big) \Bigg]  \Bigg]\]
\end{theoreme}

\begin{preuve}
	 
	\begin{align*}
	|L(z,\chi)|^{2}
	& =  \sum_{n=1}^{\infty} \sum_{q=1}^{\infty} \frac{\chi(n)}{n^{z}}
	\overline{\frac{\chi(q)}{q^{z}}} \\
	& = \sum_{\underset{m \wedge k =1}{m=1 k=1}} 
	\sum_{{p=1}}^{+\infty} \frac{\chi(mp)}{(mp)^{z}}
	\overline{\frac{\chi(kp)}{(kp)^{z}}} \\
	& = \sum_{m \wedge k =1} \sum_{{p=1}}^{+\infty} \frac{\chi(p)\chi(m)}{(p)^{z}(m)^{z}} \overline{
		\frac{\chi(p)\chi(k)}{(p)^{z}(k)^{z}} } \\
	& = \sum_{m \wedge k =1} \sum_{{p=1}}^{+\infty} \frac{\chi(p)\overline{\chi(p)}}{(p)^{2x}} \Bigg[ \frac{\chi(m)}{(m)^{z}} \overline{
		\frac{\chi(k)}{(k)^{z}}} \Bigg] \\
	& = \sum_{m \wedge k =1} \sum_{{p=1}}^{+\infty} \frac{|\chi(p)|^{2}}{(p)^{2x}} \Bigg[ \frac{\chi(m)}{(m)^{z}} \overline{
		\frac{\chi(k)}{(k)^{z}}} \Bigg] \\
	& = \sum_{m \wedge k =1} \sum_{{p=1}}^{+\infty} \frac{|\chi(p)|^{2}}{(p)^{2x}} \Bigg[ \frac{\chi(m)\overline{\chi(k)}}{(mk)^{x} {\frac{m}{k}}^{iy}} \Bigg] \\
	& = \Bigg[ \sum_{p=1}^{\infty} \frac{|\chi(p)|^{2}}{p^{2x}} \Bigg] \Bigg[ \sum_{l=1}^{\infty} \frac{1}{l^{x}} \sum_{\underset{mk=l}{m \wedge k = 1} } \frac{\chi(m)\overline{\chi(k)}}{ {\frac{m}{k}}^{iy}} \Bigg] \\ 
	& = \Bigg[ \sum_{p=1}^{\infty} \frac{|\chi(p)|^{2}}{p^{2x}} \Bigg] \Bigg[ \sum_{l=1}^{\infty} \frac{1}{l^{x}} \sum_{\underset{mk=l}{m \wedge k = 1} } \frac{\chi(m)}{{m}^{iy}}\overline{\frac{\chi(k)}{{k}^{iy}}}  \Bigg] \\
	& = \Bigg[ \sum_{p=1}^{\infty} \frac{|\chi(p)|^{2}}{p^{2x}} \Bigg] \Bigg[ \sum_{l=1}^{\infty} \frac{1}{l^{x}}  \Bigg[ \delta_{l=1}(l) + \sum_{\underset{mk=l, m>k }{m \wedge k = 1} } 2\Re\Big(\frac{\chi(m)}{{m}^{iy}}\overline{\frac{\chi(k)}{{k}^{iy}}}\Big) \Bigg]  \Bigg] 
	\end{align*}
\end{preuve}

\begin{defini}
	Le produit scalaire qui suit sera utilisé :
	
	\[
	\forall x,y,w,z \in \mathbb{R} : ~ 
	\begin{bmatrix}
		x  \\
		y
	\end{bmatrix} \cdot 
	\begin{bmatrix}
		w  \\
		z
	\end{bmatrix} = x \cdot w + y \cdot z
	\]
\end{defini}

\begin{proposition}
	\textnormal{ Le produit scalaire précédant peut être écrit sous les formes suivantes : }
	\[
	\forall x,y,w,z \in \mathbb{R} : ~ 
	x \cdot w + y \cdot z =
	\det \begin{bmatrix}
	x & i \cdot z \\
	i \cdot y & w
	\end{bmatrix} = 	
	\det \begin{bmatrix}
	x & -z \\
	y & w
	\end{bmatrix} = 	\begin{bmatrix}
	x  \\
	y
	\end{bmatrix} \cdot 
	\begin{bmatrix}
	w  \\
	z
	\end{bmatrix} 
	\]
\end{proposition}

\begin{theoreme}
	\[
\forall l \in \mathbb{N^{\star}} , \forall y \in \mathbb{R}
\Bigg[ \delta_{l=1}(l) + \sum_{\underset{mk=l, m>k }{m \wedge k = 1} } 2\Re\Big(\frac{\chi(m)}{{m}^{iy}}\overline{\frac{\chi(k)}{{k}^{iy}}}\Big) \Bigg] = \prod_{p|l}2  \cdot 
\begin{bmatrix}
\cos(y\ln(p^{n}))  \\
\sin(y\ln(p^{n}))
\end{bmatrix} \cdot 
\begin{bmatrix}
\Re(\chi(p^{n}))  \\
\Im(\chi(p^{n}))
\end{bmatrix}
\]
\end{theoreme}

\begin{preuve}
	posant 
	\[ Q(l) = \Bigg[ \delta_{l=1}(l) + \sum_{\underset{mk=l, m>k }{m \wedge k = 1} } 2\Re\Big(\frac{\chi(m)}{{m}^{iy}}\overline{\frac{\chi(k)}{{k}^{iy}}}\Big) \Bigg]\]
	En développant la partie réelle  $ \Re\Big(\frac{\chi(m)}{{m}^{iy}}\overline{\frac{\chi(k)}{{k}^{iy}}}\Big) $, et en prenant compte de la condition $ m>k $ (dont le seul rôle est d'éliminer la répétition) , on peut factoriser la fonction $Q$ : 
	\begin{align*}
	\forall p \in \mathbb{P} , \forall n \in \mathbb{N^{\star}} ~ Q(lp^{n})
	& = \frac{1}{2} \cdot \sum_{\underset{mk=l}{m \wedge k = 1} } \frac{\chi(mp^{n})}{{mp^{n}}^{iy}}\overline{\frac{\chi(k)}{{k}^{iy}}} + \frac{\chi(k)}{{k}^{iy}}\overline{\frac{\chi(mp^{n})}{{mp^{n}}^{iy}}} + \frac{\chi(kp^{n})}{{kp^{n}}^{iy}}\overline{\frac{\chi(m)}{{m}^{iy}}} + \frac{\chi(m)}{{m}^{iy}}\overline{\frac{\chi(kp^{n})}{{kp^{n}}^{iy}}} \\
	& = \frac{1}{2} \cdot [\sum_{\underset{mk=l}{m \wedge k = 1} } \frac{\chi(mp^{n})}{{mp^{n}}^{iy}}\overline{\frac{\chi(k)}{{k}^{iy}}} +
	\sum_{\underset{mk=l}{m \wedge k = 1} } \frac{\chi(k)}{{k}^{iy}}\overline{\frac{\chi(mp^{n})}{{mp^{n}}^{iy}}} + \\
	&\sum_{\underset{mk=l}{m \wedge k = 1} } \frac{\chi(kp^{n})}{{kp^{n}}^{iy}}\overline{\frac{\chi(m)}{{m}^{iy}}} +
	\sum_{\underset{mk=l}{m \wedge k = 1} } \frac{\chi(m)}{{m}^{iy}}\overline{\frac{\chi(kp^{n})}{{kp^{n}}^{iy}}}] \\
	& = \frac{1}{2} \cdot [2 \cdot \sum_{\underset{mk=l}{m \wedge k = 1} } \frac{\chi(mp^{n})}{{mp^{n}}^{iy}}\overline{\frac{\chi(k)}{{k}^{iy}}} +
	2 \cdot \sum_{\underset{mk=l}{m \wedge k = 1} } \frac{\chi(k)}{{k}^{iy}}\overline{\frac{\chi(mp^{n})}{{mp^{n}}^{iy}}}] \\
	& = \sum_{\underset{mk=l}{m \wedge k = 1} } \frac{\chi(mp^{n})}{{mp^{n}}^{iy}}\overline{\frac{\chi(k)}{{k}^{iy}}} + \frac{\chi(k)}{{k}^{iy}}\overline{\frac{\chi(mp^{n})}{{mp^{n}}^{iy}}} \\
	& = \sum_{\underset{mk=l }{m \wedge k = 1} } \frac{\chi(p^{n})}{{p^{n}}^{iy}} \frac{\chi(m)}{{m}^{iy}}\overline{\frac{\chi(k)}{{k}^{iy}}} + \frac{\chi(k)}{{k}^{iy}}\overline{ \frac{ \chi(m) }{ {m}^{iy} }}  \overline{ \frac{ \chi(p^{n}) }{ {p^{n}}^{iy} } }  \\
	& = \Bigg[ \sum_{\underset{mk=l }{m \wedge k = 1} } \frac{\chi(k)}{{k}^{iy}}\overline{ \frac{ \chi(m) }{ {m}^{iy} }} \Bigg] \Bigg[	
	\frac{\chi(p^{n})}{{p^{n}}^{iy}} + \overline{ \frac{ \chi(p^{n}) }{ {p^{n}}^{iy} } } \Bigg]    \\
	& = \Bigg[ \sum_{\underset{mk=l }{m \wedge k = 1} } \frac{\chi(k)}{{k}^{iy}}\overline{ \frac{ \chi(m) }{ {m}^{iy} }} \Bigg] \Bigg[	
	2 \cdot \Re(\frac{\chi(p^{n})}{{p^{n}}^{iy}}) \Bigg]  \\
	& = \Bigg[ \sum_{\underset{mk=l }{m \wedge k = 1} } \frac{\chi(k)}{{k}^{iy}}\overline{ \frac{ \chi(m) }{ {m}^{iy} }} \Bigg] \Bigg[	
	2  \cdot 
	\begin{bmatrix}
	\cos(y\ln(p^{n}))  \\
	\sin(y\ln(p^{n}))
	\end{bmatrix} \cdot 
	\begin{bmatrix}
	\Re(\chi(p^{n}))  \\
	\Im(\chi(p^{n}))
	\end{bmatrix}
	\Bigg] 
	\end{align*}
	Par récurrence sur chaque puissance de nombre premier dans la décomposition de n en produit de nombre premier : 
	\[ Q(n) = \prod_{p|n} 2 \cdot \begin{bmatrix}
	\cos(y\ln(p^{v_{p}(n)}))  \\
	\sin(y\ln(p^{v_{p}(n)}))
	\end{bmatrix} \cdot \begin{bmatrix}
	\Re(\chi(p^{v_{p}(n)}))  \\
	\Im(\chi(p^{v_{p}(n)}))
	\end{bmatrix}
	\]

\end{preuve}

\begin{remark}
	Cette formulation de $ Q $ prouve qu'elle est une fonction multiplicative .
	
\end{remark}

\begin{proposition}
	\textnormal{ Il est peut être intéressant de voir qu'on peut écrire la fonction $ Q $ des manières suivantes : }
	
\begin{align*}
	Q(n)
	&= \prod_{p|n} 2 \cdot [ \cos(y\ln(p^{v_{p}(n)})) \cdot \Re(\chi(p^{v_{p}(n)})) + \sin(y\ln(p^{v_{p}(n)})) \cdot  \Im(\chi(p^{v_{p}(n)})) ] \\
	&= \prod_{p|n} 2 \cdot \det \begin{bmatrix}
	\cos(y\ln(p^{v_{p}(n)})) &  -\Im(\chi(p^{v_{p}(n)})) \\
	\sin(y\ln(p^{v_{p}(n)})) & \Re(\chi(p^{v_{p}(n)}))
	\end{bmatrix} \\
	&= \prod_{p|n}  \begin{bmatrix}
	\sqrt{2} \cdot \cos(y\ln(p^{v_{p}(n)}))  \\
	\sqrt{2} \cdot \sin(y\ln(p^{v_{p}(n)}))
	\end{bmatrix} \cdot \begin{bmatrix}
	\sqrt{2} \cdot \Re(\chi(p^{v_{p}(n)}))  \\
	\sqrt{2} \cdot \Im(\chi(p^{v_{p}(n)}))
	\end{bmatrix} \\
	&= \prod_{p|n} 2 \cdot \det \begin{bmatrix}
	\cos(y\ln(p^{v_{p}(n)})) & i \cdot \Im(\chi(p^{v_{p}(n)})) \\
	i \cdot \sin(y\ln(p^{v_{p}(n)})) & \Re(\chi(p^{v_{p}(n)}))
	\end{bmatrix} \\
	&= \prod_{p|n} \det \begin{bmatrix}
	\sqrt{2} \cdot \cos(y\ln(p^{v_{p}(n)})) & i \cdot \sqrt{2} \cdot \Im(\chi(p^{v_{p}(n)})) \\
	i \cdot \sqrt{2} \cdot \sin(y\ln(p^{v_{p}(n)})) & \sqrt{2} \cdot \Re(\chi(p^{v_{p}(n)}))
	\end{bmatrix} \\
	&= \det \Bigg[ \prod_{p|n} \begin{bmatrix}
	\sqrt{2} \cdot \cos(y\ln(p^{v_{p}(n)})) & i \cdot \sqrt{2} \cdot \Im(\chi(p^{v_{p}(n)})) \\
	i \cdot \sqrt{2} \cdot \sin(y\ln(p^{v_{p}(n)})) & \sqrt{2} \cdot \Re(\chi(p^{v_{p}(n)}))
	\end{bmatrix}  \Bigg]  
\end{align*}
	
\end{proposition}

\begin{remark}
	L'ensemble des matrices :
	\[ \begin{bmatrix}
	 	c & i b \\
	 	i s & a
	\end{bmatrix} ~ ~ a , b ,c ,s \in \mathbb{R} 
	 \]
	 est stable par multiplication , son déterminant est réel. Son déterminant se calcule à l'aide des valeurs propres comme suit :
	 \[\det\begin{bmatrix}
	 	c & i b \\
	 	i s & a
	 \end{bmatrix} = \frac{c+d+\sqrt{(c+d)^2 - 4\cdot(ca+sb)}}{2} \cdot \frac{c+d-\sqrt{(c+d)^2 - 4\cdot(ca+sb)}}{2} \]
	 Ceci ne sera pas utilisé plus tard , il est mentionné simplement comme remarque.
\end{remark}

\begin{theoreme}
	
	\[ |L(z,\chi)|^{2} = \Bigg[ \sum_{p=1}^{\infty} \frac{|\chi(p)|^{2}}{p^{2x}} \Bigg] \prod_{q \in \mathbb{P}} \Bigg[ 1 + \sum_{m=1}^{\infty} \frac{1}{{q}^{mx}} \cdot 2 \cdot  \begin{bmatrix}
	\cos(y\ln(q^{m}))  \\
	\sin(y\ln(q^{m}))
	\end{bmatrix}  \cdot                \begin{bmatrix}
	\Re(\chi(q^{m}))  \\
	\Im(\chi(q^{m}))
	\end{bmatrix}
	\Bigg]\]
\end{theoreme}

\begin{preuve}
	\begin{align*}
	\forall q \in \mathbb{P} : \\ |L(z,\chi)|^{2}
	& = \Bigg[ \sum_{p=1}^{\infty} \frac{|\chi(p)|^{2}}{p^{2x}} \Bigg] \Bigg[ \sum_{l=1}^{\infty} \frac{1}{l^{x}}  \cdot \prod_{p|l} 2 \cdot  \begin{bmatrix}
	\cos(y\ln(p^{v_{p}(l)}))  \\
	\sin(y\ln(p^{v_{p}(l)}))
	\end{bmatrix}  \cdot                \begin{bmatrix}
	\Re(\chi(p^{v_{p}(l)}))  \\
	\Im(\chi(p^{v_{p}(l)}))
	\end{bmatrix}  \Bigg] \\
	& = \Bigg[ \sum_{p=1}^{\infty} \frac{|\chi(p)|^{2}}{p^{2x}} \Bigg] \Bigg[ \sum_{l=1 ,  q \nmid l}^{\infty} \frac{1}{l^{x}}  \Bigg[ \prod_{p|l} 2 \cdot  \begin{bmatrix}
	\cos(y\ln(p^{v_{p}(l)}))  \\
	\sin(y\ln(p^{v_{p}(l)}))
	\end{bmatrix}  \cdot                \begin{bmatrix}
	\Re(\chi(p^{v_{p}(l)}))  \\
	\Im(\chi(p^{v_{p}(l)}))
	\end{bmatrix} + \\
	& \sum_{l=1 , q \nmid l}^{\infty} \sum_{m=1}^{\infty} \frac{1}{l^{x}} \frac{1}{{q}^{mx}}  \cdot \prod_{p|l} 2 \cdot  \begin{bmatrix}
	\cos(y\ln(p^{v_{p}(l)}))  \\
	\sin(y\ln(p^{v_{p}(l)}))
	\end{bmatrix}  \cdot                \begin{bmatrix}
	\Re(\chi(p^{v_{p}(l)}))  \\
	\Im(\chi(p^{v_{p}(l)}))
	\end{bmatrix}  \Bigg[  2 \cdot  \begin{bmatrix}
	\cos(y\ln(q^{m}))  \\
	\sin(y\ln(q^{m}))
	\end{bmatrix}  \cdot                \begin{bmatrix}
	\Re(\chi(q^{m}))  \\
	\Im(\chi(q^{m}))
	\end{bmatrix} \Bigg]  \Bigg] \\
	& = \Bigg[ \sum_{p=1}^{\infty} \frac{|\chi(p)|^{2}}{p^{2x}} \Bigg] \Bigg[ \sum_{l=1 ,  p \nmid l}^{\infty} \frac{1}{l^{x}}  \Bigg[ \prod_{p|l} 2 \cdot  \begin{bmatrix}
	\cos(y\ln(p^{v_{p}(l)}))  \\
	\sin(y\ln(p^{v_{p}(l)}))
	\end{bmatrix}  \cdot                \begin{bmatrix}
	\Re(\chi(p^{v_{p}(l)}))  \\
	\Im(\chi(p^{v_{p}(l)}))
	\end{bmatrix} 
	\Bigg] \\                    
	& \Bigg[1+ \sum_{m=1}^{\infty} \frac{1}{{q}^{mx}} \cdot 2 \cdot   \begin{bmatrix}
	\cos(y\ln(q^{m}))  \\
	\sin(y\ln(q^{m}))
	\end{bmatrix}  \cdot                \begin{bmatrix}
	\Re(\chi(q^{m}))  \\
	\Im(\chi(q^{m}))
	\end{bmatrix}
	\Bigg] \\
	\end{align*}
	Par récurrence sur les nombres premier , on obtient la factorisation suivante : 
	
	\[ |L(z,\chi)|^{2} = \Bigg[ \sum_{p=1}^{\infty} \frac{|\chi(p)|^{2}}{p^{2x}} \Bigg] \prod_{q \in \mathbb{P}} \Bigg[ 1 + \sum_{m=1}^{\infty} \frac{1}{{q}^{mx}} \cdot 2 \cdot  \begin{bmatrix}
	\cos(y\ln(q^{m}))  \\
	\sin(y\ln(q^{m}))
	\end{bmatrix}  \cdot                \begin{bmatrix}
	\Re(\chi(q^{m}))  \\
	\Im(\chi(q^{m}))
	\end{bmatrix}
	\Bigg]\]
	 
	Ce produit correspond (après vérification dans le cas de zeta) au produit d'Euler pour $ |L(z,\chi)|^{2} $ .
\end{preuve}

\section{Preuve dans le cas de $|\zeta|^{2}$ :}

\begin{remark}
	La formule ci dessous n'est qu'un cas particulier de la factorisation du module des séries L qui la précède , néanmoins , j'ai choisi de laisser cette preuve que j'avais trouvé avant de généraliser car c'est une version moins abstraite .
\end{remark}

\begin{theoreme}
	\[
	|\zeta(z)|^{2} = \zeta(2x)\sum_{m=1}^{+\infty}\frac{1}{m^{x}}2^{\omega(m)}\prod_{j=1}^{\omega(m)}\cos(y\ln(p_{j}^{v_{p_{j}}(m)})) 
	\]
\end{theoreme}
\begin{preuve}
	La somme sur les factorisations possibles de m en nombre premier entre eux peut s'écrire plus simplement , si on écrit la décomposition de m en nombre premier :
	\[ m=\prod_{i=1}^{\omega(m)} p_{i}^{v_{p_{i}}(m)}\]
	\[ \text{rad}(m)=\prod_{i=1}^{\omega(m)} p_{i}\]
	prendre n et q premier entre eux et de produit égal à m revient à choisir un diviseur du radical de m .
	
	$ \omega(m) ~ : ~ $ est le nombre de nombre premiers distincts dans la décomposition de m
	
	le nombre de diviseur du radical de m est : $ 2^{\omega(m)} $
	\newline
	d'ou :
	\[ \sum_{\underset{n \wedge q=1}{nq=m}}\cos(y\ln(\frac{n}{q}))) = 
	\sum_{i=1}^{2^{\omega(m)}}\cos(\sum_{\varepsilon \in S}\varepsilon_{j}y\ln(p_{j})))\]
	avec $ S = \{1,-1\}^{\omega(m)}$
	\newline
	en utilisant la formule trigonométrique suivante : 
	\[ 
	\begin{aligned}
	\prod_{k=1}^n \cos \theta_k & = \frac{1}{2^n}\sum_{e\in S} \cos(e_1\theta_1+\cdots+e_n\theta_n) \\[6pt]& \text{avec }S=\{1,-1\}^n
	\end{aligned}\] 
	on obtient :
	\newline 
	\[ \sum_{i=1}^{2^{\omega(m)}}\cos(\sum_{\varepsilon \in S}\varepsilon_{j}y\ln(p_{j})))  = 2^{\omega(m)}\prod_{j=1}^{\omega(m)} \cos (y\ln(p_{j}^{v_{p_{j}}(m)}))\]
	avec $ S = \{1,-1\}^{\omega(m)}$
	\newline
	ce qui permet de réécrire la formule d'avant :
	\begin{align*}
	|\zeta(z)|^{2}
	& = \zeta(2x)\sum_{m=1}^{+\infty}\frac{1}{m^{x}}2^{\omega(m)}\prod_{j=1}^{\omega(m)}\cos(y\ln(p_{j}^{v_{p_{j}}(m)})) \\
	\end{align*}
	
\end{preuve}

\begin{remark}
	Cette formule permet de déduire un résultat connu de Hardy dans le cas ou $ y=0 $ \cite{Hardy} :
	\[ \zeta(x)^{2} = \zeta(2x)\sum_{m=1}^{+\infty} 
	\frac{2^{\omega(m)}}{m^{x}}\]
\end{remark}

\section{Cas de $ \frac{1}{2} < \Re(s) < 1 $ :}

Un Calcul similaire peut être fait en utilisant la formule suivante :
\[ \forall s \in \mathbb{C}-{1} , \Re(s) > 0 ~ : ~ \zeta(s) = \frac{1}{1-2^{1-s}}\cdot \sum_{n=1}^{\infty} \frac{(-1)^{n-1}}{n^{s}}   \]

Avec  $ \frac{1}{2} < \Re(s) < 1 $ (entre $ 0 $ et  $ \frac{1}{2} $ c'est encore plus problématique) .

Néanmoins , cela ne marche pas notamment car on a pas la convergence absolue . Ce qui complexifie (voir peut être rend impossible ) certaines manipulations sur les séries.

\begin{remark}
	
les séries factorisées dans ce chapitre continuent souvent à donner la bonne limite numériquement au sens de Cesàro pour $ \frac{1}{2} < \Re(s) < 1 $  , même quand elles ne sont pas convergentes . des séries similaires apparaissent dans \cite{LeClair5}
	
Trouver des séries convergentes (sens commun de convergence) à partir de ces séries convergentes au sens de Cesàro est une des motivations derrière la recherche qui est faite dans le chapitre 2. 

Une des pistes explorées sans sucés est de trouver une permutation de la série qui converge vers la limite de Césaro , une autre piste était de construire des mesures dépendantes de l'ordre de sommation des séries .
\end{remark}
\begin{remark}
	
L'hypothèse de Riemann découle de la preuve que la fonction zêta ne s'annule pas dans la région $ \frac{1}{2} < \Re(s) < 1 $ grâce à l'équation fonctionnelle qui implique que l'ensemble des zéros non triviaux est symétrique selon la droite $  \Re(s) = \frac{1}{2}  $ .
	
L'idée est alors que le produit d'Euler et l'équation fonctionnelle se complètent en montrant respectivement l'absence de zéro dans $ \frac{1}{2} < \Re(s) < 1 $ puis dans $ 0 < \Re(s) < \frac{1}{2} $ . le défi est de généraliser le produit d'euler dans la moitié droite de la bande critique .

\end{remark}

\begin{remark}
	Une voie possible serais de s'intéresser au permutations qui rendent la série convergente vers la limite de Césaro. 
\end{remark}

% end of body of thesis comes here

%% file: Chapter2/chapter2.tex
\chapter{Construction de l'anneau $ ( \mathds{M} , \square , \times ) $  :}
\section{Démarche :}

Sur ce chapitre , le but est de construire un certain anneau de fonction ainsi que de prouver son unicité dans un certain sens . Les opérations utilisées dans l'anneau final ont été étudiées dans plusieurs articles dont \cite{Paugam} , \cite{SITARAMAIAH} , \cite{Eckford} , \cite{VAIDYANATHASWAMY} .

On va alors commencer par quelques définitions , puis on va chercher les conditions imposées par la définition d'anneau commutatif sur une fonction poids associer à une l'opération somme de l'anneau recherché  .

\section{Définition : }

\begin{defini}
	On définit l'ensemble des fonctions multiplicatives comme suivant :
	\[
	F \in \mathbb{M} \iff 
	\begin{array}{l rcl}
	{F}  :   \mathbb{N}^{\star}  \longrightarrow  \mathbb{C} \\
	F(1)=1 \\
	a \wedge b =1 \Rightarrow F(ab)=F(a)F(b) \\
	\end{array}
	\]
\end{defini}

\begin{defini}
	On introduit les fonctions poids   $ W $  comme suivantes :
	$$
	\begin{array}{l rcl}
	W  :   {\mathbb{N}^{\star}}^{2}  \longrightarrow  \mathbb{C} \\
	\end{array}
	$$
	
\end{defini}

\begin{defini}
	On défini l'opération $\, \underset{w}{\scalebox{0.6}{$\square$}} \, $ , qu'on appellera pour l'instant W-convolution :
	$$ \forall F,G \in \mathbb{M} , \forall m \in \mathbb{N}^{\star}  :   [F  \, \underset{w}{\scalebox{0.6}{$\square$}} \, G](m) = \sum_{ab=m} F(a)G(b)W(a,b) $$

\end{defini} 

\begin{remark}
	pour alléger la notation , on notera cette opération simplement $\, \scalebox{0.6}{$\square$} \, $ .
\end{remark}

\begin{defini}
	On définit l'opération (multiplication) $ \times $ comme suivant  : 
	$$ \forall F,G \in \mathbb{M} , \forall m \in \mathbb{N}^{\star}  :   [F \times G](m) = F(m)G(m) $$ 
	
\end{defini}

\begin{defini}
	On définit les fonctions multiplicatives indicatrices comme suit  :  
	$$ \forall S \subset \mathbb{P} \times \mathbb{N}^{\star}  :   \mathds{1}_{S}(p^n) = \left\{
	\begin{array}{ll}
	1 & \mbox{si } (p,n) \in S\cup\{1,1\} \\
	0 & \mbox{sinon.}
	\end{array}
	\right. $$ 
	on écrira aussi $ \mathds{1}_{s} $ pour $ s \in \mathbb{N} $ , qu'on définit comme suit  :  $$ \mathds{1}_{s}=\mathds{1}_{S} \iff s=\prod_{(p,n)\in S} p^{n}$$
\end{defini}
\section{Construction de l'anneau  $( \mathds{M},  \square  , \times )$  : }

\begin{lemma}
	\textnormal{L'opération $ \, \scalebox{0.6}{$\square$} \, $ est commutative si et seulement si la fonction $ W(.,.) $ est commutative : }
	$$ \forall F,G \in \mathbb{M} , \forall m \in \mathbb{N}^{\star}  :   F \, \scalebox{0.6}{$\square$} \, G(m) = G \, \scalebox{0.6}{$\square$} \, F(m)   \iff \forall a,b \in \mathbb{N}^{\star}  :   W(a,b) = W(b,a) $$
\end{lemma}

\begin{preuve}
	Commençons par l'implication direct  $\Longrightarrow$  : 
	\begin{align*}
	\mathds{1}_{a}  \, \scalebox{0.6}{$\square$} \,  \mathds{1}_{b}(ab) &= \sum_{nq=ab}\mathds{1}_{a}(n)\mathds{1}_{b}(q)W(n,q) \\
	&= \mathds{1}_{a}(a)\mathds{1}_{b}(b)W(a,b) \\
	&= W(a,b) 
	\end{align*}
	Et dans l'autre sens  : 
	\begin{align*}
	\mathds{1}_{b} \, \scalebox{0.6}{$\square$} \, \mathds{1}_{a}(ab)
	&= \sum_{nq=ab}\mathds{1}_{b}(n)\mathds{1}_{a}(q)W(n,q)\\
	&= \mathds{1}_{b}(b)\mathds{1}_{a}(a)W(b,a) \\
	&= W(b,a)
	\end{align*}
	D'où le premier sens . 
	\newline
	établissons maintenant la réciproque $ \Longleftarrow $, on suppose :

	$$ \forall a,b \in \mathbb{N}^{\star}  :   W(a,b) = W(b,a) $$

	On a , $ \forall F,G \in $ $\mathbb{M} $ $ , \forall m \in $ $ \mathbb{N}^{\star} $ $  $  : 
	\begin{align*}
	F \, \scalebox{0.6}{$\square$} \, G(m)
	&= \sum_{nq=m}F(n)G(q)W(n,q)\\
	&= \sum_{\underset{n>q}{nq=m}}F(n)G(q)W(n,q)+\sum_{\underset{n<q}{nq=m}}F(n)G(q)W(n,q)+\sum_{\underset{n=q}{nq=m}}F(n)G(q)W(n,q) \\
	&= \sum_{\underset{q>n}{qn=m}}G(q)F(n)W(q,n)+\sum_{\underset{q<n}{qn=m}}G(q)F(n)W(q,n)+\sum_{\underset{q=n}{qn=m}}G(q)F(n)W(q,n) \\
	&= \sum_{nq=m}G(n)F(q)W(n,q)\\
	&= G \, \scalebox{0.6}{$\square$} \, F(m)
	\end{align*}
	D'où la réciproque . \\
	ce qui finit la preuve .
\end{preuve}

\begin{lemma}
	\textnormal{$ (\mathbb{M}, \, \scalebox{0.6}{$\square$} \,) $ est stable si et seulement si la fonction W est multiplication à deux variables selon le sens suivant  : } 
	\[ \forall F,G \in \mathbb{M}  :   F \, \scalebox{0.6}{$\square$} \, G \in \mathbb{M} \iff \forall a,b,c,d \in \mathbb{N}^{\star} ,\, ab \, \wedge \, cd = 1 :  W(a,b)W(c,d)=W(ac,bd) \]
\end{lemma}

\begin{preuve}
	Commençons par l'implication direct  $\Longrightarrow$  : 
	
	$ \forall a_{1},a_{2},b_{1},b_{2} \in \mathbb{N^{\star}} $ posant  :   $  a_{1}a_{2}=a $ et $ b_{1}b_{2}=b $ .
	
	Supposons $ a \wedge b = 1 $.
	\begin{align*}
	\mathds{1}_{a_{1}b_{1}} \, \scalebox{0.6}{$\square$} \, \mathds{1}_{a_{2}b_{2}}(ab)
	&= \sum_{nq=ab}\mathds{1}_{a}(n)\mathds{1}_{b}(q)W(n,q)\\
	&= \mathds{1}_{a}(a_{1}a_{2})\mathds{1}_{b}(b_{1}b_{2})W(a_{1}a_{2},b_{1}b_{2}) \\
	&=W(a_{1}a_{2},b_{1}b_{2})
	\end{align*}
	Vu que $ \mathds{1}_{a_{1}b_{1}} \, \scalebox{0.6}{$\square$} \, \mathds{1}_{a_{2}b_{2}} $ est une fonction multiplicative  : 
	\begin{align*}
	\mathds{1}_{a_{1}b_{1}} \, \scalebox{0.6}{$\square$} \, \mathds{1}_{a_{2}b_{2}}(ab) 
	&= [\mathds{1}_{a_{1}b_{1}} \, \scalebox{0.6}{$\square$} \, \mathds{1}_{a_{2}b_{2}}](a)\times [\mathds{1}_{a_{1}b_{1}} \, \scalebox{0.6}{$\square$} \, \mathds{1}_{a_{2}b_{2}}](b) \\
	&= [\sum_{nq=a}\mathds{1}_{a_{1}b_{1}}(n)\mathds{1}_{a_{2}b_{2}}(q)W(n,q)]\times [\sum_{nq=b}\mathds{1}_{a_{1}b_{1}}(n)\mathds{1}_{a_{2}b_{2}}(q)W(n,q)] \\
	&= [\mathds{1}_{a_{1}b_{1}}(a_{1})\mathds{1}_{a_{2}b_{2}}(a_{2})W(a_{1},a_{2})]\times[\mathds{1}_{a_{1}b_{1}}(b_{1})\mathds{1}_{a_{2}b_{2}}(b_{2})W(b_{1},b_{2})] \\
	&= W(a_{1},a_{2})W(b_{1},b_{2})
	\end{align*}
	\textnormal{ D'où le premier sens : }
	\[ W(a_{1}a_{2},b_{1}b_{2}) = W(a_{1},a_{2})W(b_{1},b_{2}) \]
	établissons maintenant la réciproque $ \Longleftarrow $  : 
	\\
	\textnormal{On suppose : } 
	\[ \forall a,b,c,d \in \mathbb{N}^{\star} ,\, ab \, \wedge \, cd = 1  :   W(a,b)W(c,d)=W(ac,bd) \]
	\textnormal{On a : }
	\begin{align*}
	\forall F,G \in \mathbb{M} , \forall a,b \in \mathbb{N}^{\star} , ~ a \wedge b =1 \textnormal{ : } [F \, \scalebox{0.6}{$\square$} \, G](ab)
	& = \sum_{nq=ab}F(n)G(q)W(n,q)\\
	& = \sum_{nq=a_{1}a_{2}b_{1}b_{2}}F(n)G(q)W(n,q)
	\end{align*}
	\textnormal{Et dans l'autre sens : }
	\begin{align*}
	[F \, \scalebox{0.6}{$\square$} \, G](a)\times [F \, \scalebox{0.6}{$\square$} \, G](b)
	&= [\sum_{n_{1}n_{2}=a}F(n_{1})G(n_{2})W(n_{1},n_{2})]\times [\sum_{q_{1}q_{2}=b}F(q_{1})G(q_{2})W(q_{1},q_{2})]\\
	&= \sum_{n_{1}n_{2}=a}F(n_{1})G(n_{2})W(n_{1},n_{2})\sum_{q_{1}q_{2}=b}F(q_{1})G(q_{2})W(q_{1},q_{2})\\
	&= \sum_{n_{1}n_{2}=a}\sum_{q_{1}q_{2}=b}F(n_{1})G(n_{2})W(n_{1},n_{2})F(q_{1})G(q_{2})W(q_{1},q_{2})\\
	&= \sum_{n_{1}n_{2}=a}\sum_{q_{1}q_{2}=b}F(n_{1}q_{1})G(n_{2}q_{2})W(n_{1}q_{1},n_{2}q_{2})
	\end{align*}
	il reste donc à prouver  :

	\[ \sum_{nq=a_{1}a_{2}b_{1}b_{2}}F(n)G(q)W(n,q) = \sum_{n_{1}n_{2}=a}\sum_{q_{1}q_{2}=b}F(n_{1}q_{1})G(n_{2}q_{2})W(n_{1}q_{1},n_{2}q_{2}) \]

	cela revient à prouver , sachant que $ a_{1}a_{2}=a $ , $ b_{1}b_{2}=b $ et $ a \wedge b = 1 $ que le changement de variable précédent est bijectif . \\
	$
	 \textnormal{étant donné }   n_{1} , n_{2} , q_{1} , q_{2} \in \mathbb{N}^{\star}  \textnormal{ qui vérifie }  n_{1}n_{2}=a  \textnormal{ et }  q_{1}q_{2}=b  .  \\ 	$
	
	 Posons $ n= n_{1}q_{1} $ et  $ q=n_{2}q_{2} $ , on a l'égalité suivante  :

	 $$ F(n_{1}q_{1}) G(n_{2}q_{2}) W(n_{1}q_{1},n_{2}q_{2}) = F(n)G(q)W(n,q) $$
	 
	où nq=ab . \\
	 inversement , étant donné $ n, q \in \mathbb{N}^{\star} $ qui vérifie :  nq=ab. \\
	on a l'égalité suivante  : 
	\[ F(n)G(q)W(n,q)=F([n \wedge a].[n \wedge b])G([q \wedge a].[q \wedge b])W([n \wedge a].[n \wedge b],[q \wedge a].[q \wedge b]) \]
	où $ [n \wedge a].[q \wedge a]=a $ et $ [n \wedge b].[q \wedge b]=b $ . \\
	\textnormal{ce qui établit l'égalité , d'où la réciproque du théorème . }
	
\end{preuve}

\begin{lemma}
	\textnormal{en supposant la stabilité et la commutativité de la structure $ (\mathbb{M},\, \scalebox{0.6}{$\square$} \,) $ , l'élément neutre est la fonction  $ \delta_{1} = \mathds{1}_{\varnothing} $ : }
	$$ \exists E \in \mathbb{M} \, \forall F \in \mathbb{M}   :   F \, \scalebox{0.6}{$\square$} \, E = E  \, \scalebox{0.6}{$\square$} \, F = F \iff \forall (n,p) \in \mathbb{N}\times \mathbb{P} ,\, W(1,p^{n})=1 , E = \delta_{1}$$
\end{lemma}

\begin{preuve}
	$ $\newline
	par commutativité , on a  : 
	$$ \exists E \in \mathbb{M} \, \forall F \in \mathbb{M}   :   F \, \scalebox{0.6}{$\square$} \, E = E  \, \scalebox{0.6}{$\square$} \, F = F \iff \exists E \in \mathbb{M} \, \forall F \in \mathbb{M}   :   F \, \scalebox{0.6}{$\square$} \, E = F $$
	par stabilité , on se restreint à  : 
	$$ \exists E \in \mathbb{M} \, \forall F \in \mathbb{M}   :   F \, \scalebox{0.6}{$\square$} \, E = F \iff \forall F \in \mathbb{M} , \forall (n,p) \in \mathbb{N}^{\star}\times \mathbb{P} ,\,   :   [F \, \scalebox{0.6}{$\square$} \, E](p^{n}) = F(p^{n}) $$
	on calcule alors la W-convolution suivante  : 
	$ \forall (n,p) \in \mathbb{N}\times \mathbb{P}  :  $ 
	\begin{align*}
	[F \, \scalebox{0.6}{$\square$} \, E](p^{n})
	&= \sum_{nq=p^{n}}F(n)E(q)W(n,q)\\
	&= \sum_{l=0}^{l=n}F(p^{l})E(p^{n-l})W(p^{l},p^{n-l}) 
	\end{align*}
	fixons  $ p \in \mathbb{P} $ prouvons le sens directe par récurrence (forte)  : 
	$$ \forall F \in \mathbb{M} , \forall n \in \mathbb{N} ,\,   :   [F \, \scalebox{0.6}{$\square$} \, E](p^{n}) = F(p^{n}) \Longrightarrow  \forall n \in \mathbb{N}^{\star} ,\, W(1,p^{n})=1 ~ \textnormal{et} ~ E = \delta_{1} $$
	cas initial  :  $ n=0 $
	\begin{align*}
	[F \, \scalebox{0.6}{$\square$} \, E](1)
	&= W(1,1) 
	\end{align*}
	\textnormal{d'où} 
	$$ W(1,1) = 1 ~ \textnormal{ et } ~ E(1) = 1 $$
	\textnormal{ l'hypothèse de récurrence : }
	$$ \forall F \in \mathbb{M} , \exists n_{0} \in \mathbb{N} ,\,   :   [F \, \scalebox{0.6}{$\square$} \, E](p^{n_{0}}) = F(p^{n_{0}}) \Longrightarrow  \forall n \in \mathbb{N}^{\star} ,\, W(1,p^{n_{0}})=1 ~ \textnormal{et} ~ E = \delta_{1} $$
	cas n+1  : 
	
	\begin{align*}
	[F \, \scalebox{0.6}{$\square$} \, E](p^{n+1})
	&= \sum_{l=0}^{l=n+1}F(p^{l})E(p^{n+1-l})W(p^{l},p^{n+1-l}) \\
	&= F(p^{n+1})E(1)W(p^{n+1},1)+F(1)E(p^{n+1})W(1,p^{n+1}) \\
	&= W(p^{n+1},1)(F(p^{n+1})+E(p^{n+1}))
	\end{align*}
	\textnormal{donc  : }
	\[W(p^{n+1},1)(F(p^{n+1})+E(p^{n+1}))-F(p^{n+1})  
	= F(p^{n+1})(W(p^{n+1},1)-1)+E(p^{n+1})W(p^{n+1},1)\]
	comme ni l'élément neutre $ E $ , ni la fonction poids $ W $ ne dépend de $ F $ , cette expression peut être vue comme un polynôme en $ F(p^{n+1}) $ .\\
	\textnormal{Posons $ A=W(p^{n+1},1)-1 ~ B=E(p^{n+1})W(p^{n+1},1) $ , on a alors  : }
	\[ A.F(p^{n+1})+B = 0 \iff A = 0 \, \text{et} \, B=0 \]
	donc :
	\[ W(p^{n+1},1) = 1 ~ et ~ E(p^{n+1}) = 0 \]
	d'où le sens direct. \\
	\textnormal{Réciproquement : }
	\begin{align*}
	[F \, \scalebox{0.6}{$\square$} \, \delta_{1}](p^{n})
	&= \sum_{l=0}^{l=n}F(p^{l})E(p^{n-l})W(p^{l},p^{n-l}) \\
	&= F(p^{n})E(1)W(p^{n+1},1) \\
	&= F(p^{n})
	\end{align*}
	ce qui finit la preuve
\end{preuve}	

\begin{lemma}
	\textnormal{ Hypothèses : } 
	\textnormal{ $ (\mathbb{M}, \scalebox{0.6}{$\square$} ) $ est stable et commutative . } \\
	\textnormal{ l'associative de la structure $ (\mathbb{M}, \scalebox{0.6}{$\square$} ) $ est équivalente à : }
	\[  \forall F,G,H \in \mathbb{M}   :   [[F \, \scalebox{0.6}{$\square$} \, G] \, \scalebox{0.6}{$\square$} \, H] = [F \, \scalebox{0.6}{$\square$} \, [G \, \scalebox{0.6}{$\square$} \, H]] \iff \forall a,b,c \in \mathbb{N}^{\star} ,\, W(a,b)W(ab,c)=W(b,c)W(bc,a) \]
\end{lemma}

\begin{preuve}
	\textnormal{ le sens directe $ \Longrightarrow $  se prouve facilement avec les fonctions indicatrices : }
	\begin{align*}
	\forall a,b,c \in \mathbb{N}^{\star} \,\mathds{1}_{a} \, \scalebox{0.6}{$\square$} \, [\mathds{1}_{b}\, \scalebox{0.6}{$\square$} \,\mathds{1}_{c}](abc) 
	&= \sum_{efg=abc}\mathds{1}_{a}(f)\mathds{1}_{b}(g)\mathds{1}_{c}(e)W(g,e)W(ge,f) \\
	&= W(b,c)W(bc,a)
	\end{align*}
	et 
	\begin{align*}
	\forall a,b,c \in \mathbb{N}^{\star} \, [\mathds{1}_{a} \, \scalebox{0.6}{$\square$} \, \mathds{1}_{b}]\, \scalebox{0.6}{$\square$} \,\mathds{1}_{c}(abc)
	&= \sum_{efg=abc}\mathds{1}_{a}(f)\mathds{1}_{b}(g)\mathds{1}_{c}(e)W(f,g)W(fg,e) \\
	&= W(a,b)W(ab,c)
	\end{align*}
	d'où  : 
	\[ \forall a,b,c \in \mathbb{N}^{\star}  :   W(a,b)W(ab,c) = W(b,c)W(bc,a) \]
	la réciproque se prouve $ \Longleftarrow $  : 
	
	$ \forall F,G,H \in \mathbb{M} , \forall (n,p) \in \mathbb{N^{\star}}\times \mathbb{P}  :  $
	\begin{align*}
	[F \, \scalebox{0.6}{$\square$} \, G]\, \scalebox{0.6}{$\square$} \,H(p^{n})
	&= \sum_{ab=p^{n}}[F \, \scalebox{0.6}{$\square$} \, G](a)H(b)W(a,b) \\
	&= \sum_{ab=p^{n}}\sum_{cd=a}F(c)G(d)W(c,d)H(b)W(cd,b) \\
	&= \sum_{bcd=p^{n}}F(c)G(d)H(b)W(c,d)W(cd,b) 
	\end{align*}
	et  : 
	\begin{align*}
	F \, \scalebox{0.6}{$\square$} \, [G\, \scalebox{0.6}{$\square$} \,H](p^{n})
	&= \sum_{ba=p^{n}}F(b)[G \, \scalebox{0.6}{$\square$} \, H](a)W(b,a) \\
	&= \sum_{ba=p^{n}}F(b)\sum_{cd=a}G(c)H(d)W(c,d)W(b,a) \\
	&= \sum_{bcd=p^{n}}F(b)G(c)H(d)W(c,d)W(b,cd) \\
	&= \sum_{bcd=p^{n}}F(c)G(d)H(b)W(d,b)W(c,db) \\ 
	&= \sum_{bcd=p^{n}}F(c)G(d)H(b)W(d,b)W(db,c) 
	\end{align*}
	or on a  :  
	\[ \forall a,b,c \in \mathbb{N}^{\star}  :   W(a,b)W(ab,c) = W(b,c)W(bc,a) \]
	d'où : 
	$$\sum_{bcd=p^{n}}F(c)G(d)H(b)W(c,d)W(cd,b)=\sum_{bcd=p^{n}}F(c)G(d)H(b)W(d,b)W(db,c) $$
	d'ou le sens direct .
	
\end{preuve} 

\begin{lemma}
	\textnormal{ notons : } $\mathds{1}  : = \mathds{1}_{\mathbb{N}\times\mathbb{P}}   $.
	
	\textnormal{$ (\mathbb{M},\, \times \,) $ est stable , commutatif , associatif et admet $ \mathds{1} $ comme élément neutre . }
\end{lemma}

\begin{preuve}
	Stabilité : 
	$$ \forall F,G \in \mathbb{M}, \forall a,b \in \mathbb{N^{\star}}^{2}  :  [F \times G](ab)=F(ab)\times G(ab)=F(a)G(a)F(b)G(b)=[F\times G](a)\times[F\times G](b) $$
	d'où  : 
	$$ \forall F,G \in \mathbb{M}  :  F\times G \in \mathbb{M} $$
	Commutativité :  
	$$ \forall F,G \in \mathbb{M}, \forall m \in \mathbb{N^{\star}} :  [F \times G](m)=F(m)\times G(m)=G(m)\times F(m)=[G \times F](m) $$
	d'où : 
	$$ \forall F,G \in \mathbb{M}  :  F \times G = G \times F  $$
	Associativité : 
	$$ \forall F,G,H \in \mathbb{M}, \forall m \in \mathbb{N^{\star}}  :  [[F \times G] \times H ](m)=[F\times G](m) \times H(m)=F(m)G(m)H(m)=[F\times [G\times H ]](m) $$
	d'où : 
	$$ \forall F,G,H \in \mathbb{M}  :  [[F \times G] \times H = [F\times [G\times H ]]  $$
	$ \mathds{1} $ est l'élément neutre : 
	$$ \forall F \in \mathbb{M} , \forall (n,p) \in \mathbb{N}\times \mathbb{P}  :  [F\times \mathds{1}](p^{n}) = F(p^{n})\times\mathds{1}(p^{n}) = F(p^{n}) $$
	d'où ,  par unicité de l'élément neutre , $ \mathds{1} $ est l'élément neutre de  $ (\mathbb{M},\, \times \,) $ . \\
	ce qui prouve le théorème.
	
\end{preuve}

\begin{lemma}
	\textnormal{Hypothèse : }
	\[ \forall a,b,c,d \in \mathbb{N}^{\star} ,\, ab \, \wedge \, cd = 1  :   W(a,b)W(c,d)=W(ac,bd) \]
	\textnormal{Qui peut être généralisé en  : }
	\begin{align*}
	\forall n \in \mathbb{N} \, , \, \forall \, {a}_{0},{a}_{1}, ... , {a}_{n} \in \mathbb{N}^{\star} , \, &\forall \, {b}_{0},{b}_{1}, ... , {b}_{n} \in \mathbb{N}^{\star} \, | \, \forall i,j \in \llbracket 0,n\rrbracket \, , \,{a}_{i}{b}_{i} \wedge  {a}_{j}{b}_{j} = 1 \\
	& \prod_{i=0}^{i=n}W({a}_{i},{b}_{i}) = W(\prod_{i=0}^{i=n}{a}_{i},\prod_{i=0}^{i=n}{b}_{i})
	\end{align*}
\end{lemma}
\begin{preuve}
	On va procéder par récurrence : \\
	\textnormal{Cas $n=0$ : }
	\[ \prod_{i=0}^{i=0}W({a}_{i},{b}_{i}) = W(\prod_{i=0}^{i=0}{a}_{i},\prod_{i=0}^{i=0}{b}_{i}) 
	 \iff 1 = W(1,1) \]
	hypothèse de récurrence  : 
	\begin{align*}
	\exists {n}_{0} \in \mathbb{N} \, , \, \forall \, {a}_{0},{a}_{1}, ... , {a}_{{n}_{0}} \in \mathbb{N}^{\star} , \, &\forall \, {b}_{0},{b}_{1}, ... , {b}_{{n}_{0}} \in \mathbb{N}^{\star} \, | \, \forall i,j \in \llbracket 0,{n}_{0}\rrbracket \, , \,{a}_{i}{b}_{i} \wedge  {a}_{j}{b}_{j} = 1 \\
	& \prod_{i=0}^{i={n}_{0}}W({a}_{i},{b}_{i}) = W(\prod_{i=0}^{i={n}_{0}}{a}_{i},\prod_{i=0}^{i={n}_{0}}{b}_{i})
	\end{align*}
	Le cas $n+1$  : 
	\begin{align*}
	\prod_{i=0}^{i=n+1}W({a}_{i},{b}_{i}) &= [\prod_{i=0}^{i=n}W({a}_{i},{b}_{i})]W({a}_{n+1},{b}_{n+1}) \\
	&=W(\prod_{i=0}^{i=n}{a}_{i},\prod_{i=0}^{i=n}{b}_{i})W({a}_{n+1},{b}_{n+1}) \tag*{[hypothèse récurrence]}\\
	&=W(\prod_{i=0}^{i=n+1}{a}_{i},\prod_{i=0}^{i=n+1}{b}_{i})\tag*{[${a}_{n+1}{b}_{n+1} \wedge \prod_{i=0}^{i=n}{a}_{i}\prod_{i=0}^{i=n}{b}_{i}=1$]}
	\end{align*}
	\textnormal{D'où le résultat recherché. }
\end{preuve}

\begin{lemma}
	\textnormal{Hypothèse  :  $ (\mathbb{M}, \, \scalebox{0.6}{$\square$} \,) $ est stable , commutatif et il admet un élément neutre. }\\
	\textnormal{La fonction W peut s'écrire comme tel  : } 
	\begin{align*}
	\forall n,q \in \mathbb{N}^{\star}  :  W(n,q) &=\prod_{\stackrel{p | \frac{nq}{(n\wedge q)^{2}}}{p \in \mathbb{P}}} W({p}^{v_{p}(\frac{nq}{(n\wedge q)^{2}})},1)\prod_{\stackrel{p | n \wedge q}{p \in \mathbb{P}}}W({p}^{v_{p}(n \wedge q)},{p}^{v_{p}(n \wedge q)})\\
	&=\prod_{\stackrel{p | n \wedge q}{p \in \mathbb{P}}}W({p}^{v_{p}(n \wedge q)},{p}^{v_{p}(n \wedge q)})
	\end{align*}
\end{lemma}

\begin{remark}
	Cela signifie que dans ce cas , la fonction W est entièrement déterminée par l'image des couples égaux de puissance de nombre premier.
\end{remark}

\begin{preuve}
	\textnormal{la fonction W vérifie , par hypothèse , la propriété suivante : }
	\[ \forall a,b,c,d \in \mathbb{N}^{\star} ,\, ab \, \wedge \, cd = 1  :   W(a,b)W(c,d)=W(ac,bd) \]
	On a $ \frac{nq}{(n \wedge q)^{2}} \wedge (n \wedge q)^{2} = 1 $ d'où  : 
	\begin{align*}
	\forall n,q \in \mathbb{N}^{\star}  :  W(n,q)
	&=W(\frac{n}{n \wedge q},\frac{q}{n \wedge q})W(n \wedge q,n \wedge q)
	\end{align*}
	Calculons la première partie  : 
	\begin{align*}
	\forall n,q \in \mathbb{N}^{\star}  :  W(\frac{n}{n \wedge q},\frac{q}{n \wedge q}) 
	&=W(\frac{n}{n \wedge q},1)W(1,\frac{q}{n \wedge q})\\
	&=W(\frac{nq}{{(n \wedge q)}^{2}},1)\tag*{$\frac{n}{n \wedge q}\times 1 \wedge \frac{q}{n \wedge q}\times 1=1 $}\\
	&=\prod_{\underset{p \in \mathbb{P}}{p | \frac{nq}{{(n \wedge q)}^{2}}}}W({p}^{v_{p}(\frac{nq}{{(n \wedge q)}^{2}})},1)
	\end{align*}
	Calculons ensuite la seconde partie  : 
	\begin{align*}
	\forall n,q \in \mathbb{N}^{\star}  :  W(n \wedge q,n \wedge q)
	&= \prod_{\underset{p \in \mathbb{P}}{p | n \wedge q}}W({p}^{v_{p}(n \wedge q)},{p}^{v_{p}(n \wedge q)})
	\end{align*}
	d'où  : 
	\begin{align*}
	\forall n,q \in \mathbb{N}^{\star}  :  W(n,q)
	&=\prod_{\underset{p \in \mathbb{P}}{p | \frac{nq}{{(n \wedge q)}^{2}}}}W({p}^{v_{p}(\frac{nq}{{(n \wedge q)}^{2}})},1)\prod_{\underset{p \in \mathbb{P}}{p | n \wedge q}}W({p}^{v_{p}(n \wedge q)},{p}^{v_{p}(n \wedge q)})
	\end{align*}
	par hypothèse (existence de l'élément neutre) , on obtient que  : 
	\begin{align*}
	\forall n,q \in \mathbb{N}^{\star}  :  W(n,q)
	&=\prod_{\underset{p \in \mathbb{P}}{p | n \wedge q}}W({p}^{v_{p}(n \wedge q)},{p}^{v_{p}(n \wedge q)})
	\end{align*}
\end{preuve}

\begin{lemma}
	\textnormal{ Hypothèses : Commutativité , stabilité , associativité des deux opérations dans la structure $(\mathbb{M},  \scalebox{0.6}{$\square$}  , \times ) $ . } \\
	\textnormal{Dans la structure $(\mathbb{M},  \scalebox{0.6}{$\square$}  , \times ) $ , l'opération \scalebox{0.6}{$\square$}  est distributive par rapport à $ \times $  si et seulement si la fonction W , poids de la convolution , vérifie certaine condition : }
	\[ \forall F,G,H \in \mathbb{M}  :  [F \, \scalebox{0.6}{$\square$} \, G] \times H = [F \times H ] \, \scalebox{0.6}{$\square$} \, [G \times H] \iff  \forall a,b \in \mathbb{N}^{\star}  :  W(a,b) =\begin{cases}
	\begin{array}{ll}
	1 \,\textrm{ \textnormal{si} } \, a \wedge b = 1 \\
	0 \,\textrm{ \textnormal{si non}}
	\end{array}
	\end{cases}  \]
	
\end{lemma}

\begin{preuve}
	Commençons par le sens direct $ \Longrightarrow $ : 
	\begin{align*}
	\forall \, l,f  \in \mathbb{N} \, l+f=n , \forall (n,p) \in \mathbb{N^{\star}}\times \mathbb{P}  :  \\ [\mathds{1}_{{p}^{l}} \, \scalebox{0.6}{$\square$} \, \mathds{1}_{{p}^{f}}] \times \mathds{1}_{{p}^{n}} ({p}^{n}) 
	&= [\mathds{1}_{{p}^{l}} \, \scalebox{0.6}{$\square$} \, \mathds{1}_{{p}^{f}}] ({p}^{n}) \times \mathds{1}_{{p}^{n}} ({p}^{n}) \\ 
	&=\sum_{ab=p^{n}}\mathds{1}_{{p}^{l}}(a)\mathds{1}_{{p}^{f}}(b)W(a,b) \\
	&=W({p}^{l},{p}^{f})
	\end{align*}
	\begin{align*}
	\forall \, l,f  \in \mathbb{N} \, l+f=n , \forall (n,p) \in \mathbb{N^{\star}}\times \mathbb{P}  :  \\ [\mathds{1}_{{p}^{l}} \times \mathds{1}_{{p}^{n}} ] \, \scalebox{0.6}{$\square$} \, [\mathds{1}_{{p}^{f}} \times \mathds{1}_{{p}^{n}} ] ({p}^{n}) 
	&=\sum_{ab=p^{n}}[\mathds{1}_{{p}^{l}} \times \mathds{1}_{{p}^{n}} ](a) [\mathds{1}_{{p}^{f}} \times \mathds{1}_{{p}^{n}} ] (b) W(a,b) \\ 
	&=\sum_{ab=p^{n}} \mathds{1}_{{p}^{l}}(a) \mathds{1}_{{p}^{n}}(a) \mathds{1}_{{p}^{f}}(b) \mathds{1}_{{p}^{n}}(b) W(a,b) \\
	&=\mathds{1}_{{p}^{l}}({p}^{n}) \mathds{1}_{{p}^{n}}({p}^{n}) \mathds{1}_{{p}^{f}}(1) \mathds{1}_{{p}^{n}}(1) W({p}^{n},1) \\
	&+ \mathds{1}_{{p}^{l}}(1) \mathds{1}_{{p}^{n}}(1) \mathds{1}_{{p}^{f}}({p}^{n}) \mathds{1}_{{p}^{n}}({p}^{n}) W(1,{p}^{n})\\
	&=[\mathds{1}_{{p}^{l}}({p}^{n}) +  \mathds{1}_{{p}^{f}}({p}^{n})] W({p}^{n},1)
	\end{align*} 
	Cas $l.f=0$ :  \\
	on a alors $  l=n ~ \textrm{ou} ~ f=n$ : 
	\[
	[\mathds{1}_{{p}^{l}}({p}^{n}) +  \mathds{1}_{{p}^{f}}({p}^{n})] W({p}^{n},1)
	=[\mathds{1}_{{p}^{n}}({p}^{n}) +  \mathds{1}_{{p}^{0}}({p}^{n})] W({p}^{n},1) 
	=1
	\] 
	Cas $l.f \neq 0$ : \\
	on a alors  $ \, l \neq 0 \, \textrm{,} \, f \neq 0 \, \textrm{,} \,  l \neq n \, \textrm{et} \, f \neq n \, $ et :
	\begin{align*}
	[\mathds{1}_{{p}^{l}}({p}^{n}) +  \mathds{1}_{{p}^{f}}({p}^{n})] W({p}^{n},1)
	&=[0+0] W({p}^{n},1) \\
	&=0
	\end{align*}
	D'où , en résumant :
	\[
	W({p}^{l},{p}^{f}) =\begin{cases}
	\begin{array}{ll}
	1 \,\textrm{ si } \, l.f=0 \\
	0 \,\textrm{ si non}
	\end{array}
	\end{cases}
	\]
	En d'autre terme (équivalent) :
	\[
	W({p}^{l},{p}^{f}) =\begin{cases}
	\begin{array}{ll}
	1 \,\textrm{ si } \, {p}^{l} \wedge {p}^{f} = 1 \\
	0 \,\textrm{ si non}
	\end{array}
	\end{cases}
	\]
	D'après le résultat précédent , cela détermine entièrement la fonction , or la fonction : 
	\[
	\forall a,b \in \mathbb{N}^{\star}  :   W(a,b)
	\begin{cases}
	\begin{array}{ll}
	1 \,\textrm{ si } \, a \wedge b = 1 \\
	0 \,\textrm{ si non}
	\end{array}
	\end{cases}
	\]
	vérifie ces conditions  , D'où le sens directe.\\
	Passons maintenant au sens indirect $ \Longleftarrow $  : 
	\begin{align*}
	\forall F,G,H \in \mathbb{M} \, \forall (n,p) \in \mathbb{N^{\star}}\times \mathbb{P}  :  [F \, \scalebox{0.6}{$\square$} \, G] \times H ({p}^{n}) &= [\sum_{\underset{a \wedge b = 1}{ab={p}^{n}}} F(a)G(b)W(a,b)]H({p}^{n}) \\
	&=[F(1)G({p}^{n})W(1,{p}^{n})+F({p}^{n})G(1)W({p}^{n},1)]H({p}^{n})\\
	&=G({p}^{n})H({p}^{n})+F({p}^{n})H({p}^{n})
	\end{align*} 
	\begin{align*}
	\forall F,G,H \in \mathbb{M} \, \forall (n,p) \in \mathbb{N^{\star}}\times \mathbb{P}  :  [F \times H ] \, \scalebox{0.6}{$\square$} \, [G \times H] ({p}^{n}) &= \sum_{\underset{a \wedge b = 1}{ab={p}^{n}}} F(a)H(a)G(b)H(b)W(a,b) \\
	&=F(1)H(1)G({p}^{n})H({p}^{n})W(1,{p}^{n})+F({p}^{n})H({p}^{n})G(1)H(1)W({p}^{n},1)\\
	&=G({p}^{n})H({p}^{n})+F({p}^{n})H({p}^{n})
	\end{align*} 
	D'où le résultat recherché.
	
\end{preuve}

\begin{lemma}
	\textnormal{Hypothèse :
	Commutativité , stabilité , existence élément neutre dans la structure $(\mathbb{M},  \scalebox{0.6}{$\square$}  , \times ) $ ainsi que la distributivité de $\scalebox{0.6}{$\square$}$ par rapport à $\times$.} \\
	\textnormal{Dans la structure $ (\mathbb{M},  \scalebox{0.6}{$\square$}  , \times ) $  , pour tout élément de $\mathbb{M}$ , il existe un unique inverse par rapport à la loi $\scalebox{0.6}{$\square$}$ définie comme tel : }
	\[ \forall F \in \mathbb{M}  :  [F \, \scalebox{0.6}{$\square$} \, I_{F}]  = \delta_{1} \iff  \forall  n,p \in \mathbb{N}^{\star} \times \mathbb{P} \,  :  \, I_{F}({p}^{n})=-F({p}^{n}) \]
\end{lemma}

\begin{preuve}
	Calculons  : 
	\begin{align*}
	\forall F \in \mathbb{M} \forall  n,p \in \mathbb{N}^{\star} \times \mathbb{P}  :  [F \, \scalebox{0.6}{$\square$} \, I_{F}]({p}^{n}) &= \sum_{\underset{a \wedge b = 1}{ab = {p}^{n}}} F(a) I_{F}(b)W(a,b) \\
	&= F(1)I_{F}({p}^{n}) + F({p}^{n})I_{F}(1) \\
	&= I_{F}({p}^{n})+F({p}^{n}) \\
	&=\delta_{1}({p}^{n}) \\
	&=0
	\end{align*}
	D'où 
	\[ \forall F \in \mathbb{M} \forall  n,p \in \mathbb{N}^{\star} \times \mathbb{P}  :  I_{F}({p}^{n}) = -F({p}^{n})\]
	Dans le cas n=0  : 
	\begin{align*}
	\forall F \in \mathbb{M}   :  [F \, \scalebox{0.6}{$\square$} \, I_{F}](1) &= \sum_{\underset{a \wedge b = 1}{ab = 1}} F(a) I_{F}(b)W(a,b) \\
	&= F(1)I_{F}(1) \\
	&=\delta_{1}(1) \\
	&=1
	\end{align*}
	D'où
	$$ I_{F}(1)=1  $$
	Résumons :
	\begin{align*}
	\forall F \in \mathbb{M} \forall  n,p \in \mathbb{N} \times \mathbb{P}  :  I_{F}({p}^{n}) =  
	\begin{cases}
	\begin{array}{ll}
	1 \,\textrm{ si } \, n = 0 \\
	-F({p}^{n}) \,\textrm{ si non }
	\end{array}
	\end{cases}
	\end{align*}
	Ce qui représente la même fonction ci-dessous :
	\begin{align*}
	\forall F \in \mathbb{M} \forall  n \in \mathbb{N}^{\star}  :  I_{F}(n) = F^{-1}(n) = (-1)^{\omega(n)} F(n)
	\end{align*}
\end{preuve}

\section{Résultats :}

Dans cette section , on présente les principaux résultats conclus à partir des lemmes précédents :

\begin{theoreme}
	\textnormal{La W - convolution a une formule produit : } \\
	\textnormal{Pour : }
	$ \begin{array}{l rcl}
		\forall a,b \in \mathbb{N}^{\star}  :  W(a,b) = 
		\begin{cases}
			\begin{array}{ll}
				1 \,\textrm{ si } \, a \wedge b = 1 \\
				0 \,\textrm{ si non} 
			\end{array}
		\end{cases}
	\end{array} $
	\textnormal{On a :}
	
	\[ \forall F,G \in \mathbb{M} \, \forall n \in \mathbb{N}^{\star}  :  [F \, \scalebox{0.6}{$\square$} \, G](n)  = \sum_{\underset{a \wedge  b= 1}{ab = n}}F(a)G(b) = \prod_{p|n}(F({p}^{v_{p}(n)}))+G({p}^{v_{p}(n)})) \]
\end{theoreme}

\begin{preuve}
L'image des puissances de nombres premiers caractérise les fonctions multiplicatives , calculons alors ces images pour la convolution recherché  : 
\begin{align*}
\forall F,G \in \mathbb{M} \, \forall n,p \in \mathbb{N}^{\star} \times \mathbb{P}  :  [F \, \scalebox{0.6}{$\square$} \, G]({p}^{n})  
&= \sum_{\underset{a \wedge  b= 1}{ab = {p}^{n}}}F(a)G(b) \\
&=\begin{cases}
\begin{array}{ll}
1 \,\textrm{ si } \, n = 0 \\
F({p}^{n})+G({p}^{n}) \,\textrm{ si non }
\end{array}
\end{cases} 
\end{align*} 
\begin{align*}
\forall F,G \in \mathbb{M} \, \forall n,p \in \mathbb{N}^{\star} \times \mathbb{P}  :  [F \, \scalebox{0.6}{$\square$} \, G]({p}^{n})  
&=\prod_{q|{p}^{n}}(F({q}^{v_{q}({p}^{n})}))+G({q}^{v_{q}({p}^{n})})) \\
&=\begin{cases}
\begin{array}{ll}
1 \,\textrm{ si } \, n = 0 \\
F({p}^{n})+G({p}^{n}) \,\textrm{ si non }
\end{array}
\end{cases} 
\end{align*} 
D'où le résultat recherché. 
\end{preuve}

\begin{theoreme}
	\[ (\mathbb{M},  \scalebox{0.6}{$\square$}  , \times ) ~ \, \textnormal{est un anneau commutatif} \iff 
	\begin{array}{l rcl}
	\forall a,b \in \mathbb{N}^{\star}  :  W(a,b) =
	\begin{cases}
	\begin{array}{ll}
	1 \,\textnormal{ si } \, a \wedge b = 1 \\
	0 \,\textnormal{ si non}
	\end{array}
	\end{cases}
	\end{array}  \]
\end{theoreme}

\begin{preuve}
	En combinant les lemmes précédents :
\begin{align*}		
	\begin{split}
	(\mathbb{M},  \scalebox{0.6}{$\square$}  , \times ) ~ \, \textnormal{est un anneau commutatif}
	& \Rightarrow  
		\begin{cases}
			\begin{array}{l rcl}
				\forall a,b \in \mathbb{N}^{\star}  :   W(a,b) = W(b,a) \\
				\forall a,b,c,d \in \mathbb{N}^{\star} ,\, ab \, \wedge \, cd = 1  :   W(a,b)W(c,d)=W(ac,bd) \\
				\forall (n,p) \in \mathbb{N}\times \mathbb{P} ,\, W(1,p^{n})=1 \\
				\forall a,b,c \in \mathbb{N}^{\star} ,\, W(a,b)W(ab,c)=W(b,c)W(bc,a) \\
				\begin{array}{l rcl}
					\forall a,b \in \mathbb{N}^{\star}  :  W(a,b) =
					\begin{cases}
						\begin{array}{ll}
							1 \,\textnormal{ si } \, a \wedge b = 1 \\
							0 \,\textnormal{ si non}
						\end{array}
					\end{cases}
				\end{array}
			\end{array}
		\end{cases}\\
	& \Rightarrow
		\begin{cases}
			\begin{array}{l rcl}
			\forall a,b \in \mathbb{N}^{\star}  :  W(a,b) =
				\begin{cases}
					\begin{array}{ll}
						1 \,\textnormal{ si } \, a \wedge b = 1 \\
						0 \,\textnormal{ si non}
					\end{array}
				\end{cases}
			\end{array} 
		\end{cases}
	\end{split}
\end{align*}
Réciproquement (la vérification est plus simple en utilisant la formule du produit) :
\begin{align*}		
\begin{split}
\begin{array}{l rcl}
\forall a,b \in \mathbb{N}^{\star}  :  W(a,b) =
\begin{cases}
\begin{array}{ll}
1 \,\textnormal{ si } \, a \wedge b = 1 \\
0 \,\textnormal{ si non}
\end{array}
\end{cases}
\end{array} 
& \Rightarrow  
(\mathbb{M},  \scalebox{0.6}{$\square$}  , \times ) ~ \, \textnormal{est un anneau commutatif}\\
\end{split}
\end{align*}
\end{preuve}

%% file: Chapter3/chapter3.tex
\chapter{Propositions associées à l'anneau $ ( \mathds{M} , \square , \times ) $ : }

Plusieurs résultats sont présentés dans cette sections , La majorité sont déduits des opérations de l'anneau construit précédemment . 

\section{Définition :}
\begin{remark}
	Dans cette section, on va utiliser l'anneau $ (\mathbb{M},  \scalebox{0.6}{$\square$}  , \times )$ avec l'opération de convolution défini par :
	\[
	\begin{array}{l rcl}
		\forall a,b \in \mathbb{N}^{\star} \textnormal{ : } W(a,b) =
		\begin{cases}
			\begin{array}{ll}
				1 \,\textrm{ si } \, a \wedge b = 1 \\
				0 \,\textrm{ si non}
			\end{array}
		\end{cases}
	\end{array}
	\]
\end{remark}

\begin{defini}
	On définit l'ensemble des fonctions complètement multiplicatives comme suivant :
	\[
	F \in \mathbb{M}_{c} \iff 
	\begin{array}{l rcl}
	F \in \mathbb{M} \\
	a,b \in \mathbb{N^{\star}}  : F(ab)=F(a)F(b) \\
	\end{array}
	\]
\end{defini}

\begin{remark}
	$ \mathbb{M}_{c} $ est stable par multiplication
\end{remark}

\section{Résultats :}

\begin{proposition}
	\textnormal{pour F et G complètement multiplicatives :}
	\[\forall F,G \in \mathbb{M}_{c} :  D(F,s) \times D(G,s) = D(F \times G,2s) \times D( F {\, \scalebox{0.5}{$\square$} \,} G ,s) \]
\end{proposition}
\begin{preuve}
	\begin{align*}
	\forall F,G \in \mathbb{M}_{c} : D(F,s) \times D(G,s) 
	& = \Bigg[ \sum_{n=1}^{+\infty} \frac{F(n)}{n^{s}} \Bigg] \times \Bigg[ \sum_{q=1}^{+\infty} \frac{G(q)}{q^{s}} \Bigg] \\
	& =  \sum_{n=1}^{\infty} \sum_{q=1}^{\infty} \frac{F(n)}{n^{s}}
	\frac{G(q)}{q^{s}} \\
	& = \sum_{\underset{m \wedge k =1}{m=1 k=1}}  \sum_{{r=1}}^{+\infty} 
	\frac{F(mr)}{(mr)^{s}}
	\frac{G(kr)}{(kr)^{s}} \\
	& = \sum_{{r=1}}^{+\infty} \frac{F(r)G(r)}{(r)^{2s}} \sum_{m \wedge k =1}  \Bigg[ \frac{F(m)}{(m)^{s}} 
	\frac{G(k)}{(k)^{s}} \Bigg] \\
	& = \sum_{{r=1}}^{+\infty} \frac{F(r)G(r)}{(r)^{2s}} \sum_{l=1}^{\infty}\sum_{\underset{mk=l}{m \wedge k =1}}  \Bigg[ \frac{F(m)}{(m)^{s}} 
	\frac{G(k)}{(k)^{s}} \Bigg] \\
	& = \Bigg[ \sum_{r=1}^{\infty} \frac{F(r)G(r)}{r^{2s}} \Bigg] \Bigg[ \sum_{l=1}^{\infty}\sum_{\underset{mk=l}{m \wedge k =1}}  \Bigg[ \frac{F(m)}{(m)^{s}} 
	\frac{G(k)}{(k)^{s}} \Bigg] \\ 
	& = \Bigg[ \sum_{r=1}^{\infty} \frac{F(r)G(r)}{r^{2s}} \Bigg] \Bigg[ \sum_{l=1}^{\infty}\frac{ \sum_{\underset{mk=l}{m \wedge k =1}}  F(m) G(k) }{l^{s}} \Bigg] \\ 
	& = \Bigg[ \sum_{r=1}^{\infty} \frac{F(r)G(r)}{r^{2s}} \Bigg] \Bigg[ \sum_{l=1}^{\infty}\frac{ F {\, \scalebox{0.5}{$\square$} \,} G (l) }{l^{s}} \Bigg] \\ 
	& = D(F \times G,2s) \times D( F {\, \scalebox{0.5}{$\square$} \,} G ,s)
	\end{align*}
\end{preuve}

\begin{proposition}
	\[  \forall F , G \in \mathbb{M}_{c} : ~ D(F,s) \times D(G,\overline{s}) = D(F \times G,2x) \times D( \frac{F}{\text{Id}_{e}^{iy}} \scalebox{0.6}{$\square$} \frac{G}{\text{Id}_{e}^{-iy}} , x ) \]
\end{proposition}
\begin{preuve}
	\begin{align*}
	F,G \in \mathbb{M}_{c} ~ \forall y \in \mathbb{C} :  D(F,s) \times D(G,\overline{s}) 
	& = \Bigg[ \sum_{n=1}^{+\infty} \frac{F(n)}{n^{s}} \Bigg] \times \Bigg[ \sum_{q=1}^{+\infty} \frac{G(q)}{q^{\overline{s}}} \Bigg] \\
	& =  \sum_{n=1}^{\infty} \sum_{q=1}^{\infty} \frac{F(n)}{n^{s}}
	\frac{G(q)}{q^{\overline{s}}} \\
	& = \sum_{\underset{m \wedge k =1}{m=1 k=1}}  \sum_{{p=1}}^{+\infty} 
	\frac{F(mp)}{(mp)^{s}}
	\frac{G(kp)}{(kp)^{\overline{s}}} \\
	& = \sum_{\underset{m \wedge k =1}{m=1 k=1}}  \sum_{{p=1}}^{+\infty} 
	\frac{F(mp)}{(mp)^{x+iy}}
	\frac{G(kp)}{(kp)^{x-iy}} \\
	& = \sum_{{p=1}}^{+\infty} \frac{F(p)G(p)}{(p)^{2x}} \sum_{m \wedge k =1}  \Bigg[ \frac{F(m)}{(m)^{x+iy}} 
	\frac{G(k)}{(k)^{x-iy}} \Bigg] \\
	& = \sum_{{p=1}}^{+\infty} \frac{F(p)G(p)}{(p)^{2x}} \sum_{l=1}^{\infty}\sum_{\underset{mk=l}{m \wedge k =1}}  \Bigg[ \frac{F(m)}{(m)^{x+iy}} 
	\frac{G(k)}{(k)^{x-iy}} \Bigg] \\
	& = \Bigg[ \sum_{p=1}^{\infty} \frac{F(p)G(p)}{p^{2x}} \Bigg] \Bigg[ \sum_{l=1}^{\infty}\sum_{\underset{mk=l}{m \wedge k =1}}  \Bigg[ \frac{F(m)}{(m)^{x+iy}} 
	\frac{G(k)}{(k)^{x-iy}} \Bigg] \\ 
	& = \Bigg[ \sum_{p=1}^{\infty} \frac{F(p)G(p)}{p^{2x}} \Bigg] \Bigg[ \sum_{l=1}^{\infty}\frac{ \sum_{\underset{mk=l}{m \wedge k =1}}  \frac{F(m)}{(m)^{iy}} 
	\frac{G(k)}{(k)^{-iy}} }{l^{x}} \Bigg] \\ 
	& = \Bigg[ \sum_{p=1}^{\infty} \frac{F(p)G(p)}{p^{2x}} \Bigg] \Bigg[ \sum_{l=1}^{\infty}\frac{ \frac{F}{\text{Id}_{e}^{iy}} {\, \scalebox{0.5}{$\square$} \,} \frac{G}{\text{Id}_{e}^{-iy}} (l) }{l^{x}} \Bigg] \\ 
	& = D(F \times G,2x) \times D( \frac{F}{\text{Id}_{e}^{iy}} {\, \scalebox{0.5}{$\square$} \,} \frac{G}{\text{Id}_{e}^{-iy}} ,x)
	\end{align*}
\end{preuve}

\begin{proposition}
	\textnormal{pour F et G complètement multiplicatives :}
	\[ \forall F,G \in \mathbb{M_{c}} ~: ~ F \times G = \delta_{1} \Rightarrow D(F,s) \times D(G,s) = D( F {\, \scalebox{0.5}{$\square$} \,} G ,s) \]
\end{proposition}

\begin{preuve}
	\[ \forall F,G \in \mathbb{M_{c}} ~: ~  F \times G = \delta_{1} \Rightarrow D(F \times G,s) = 1 \]
	d'où :
	\[ \forall F,G \in \mathbb{M_{c}} ~: ~  D(F,s) \times D(G,s) = D( F {\, \scalebox{0.5}{$\square$} \,} G ,s) \times D(F \times G,s) = D( F {\, \scalebox{0.5}{$\square$} \,} G ,s)  \]
\end{preuve}

\begin{proposition}
	\textnormal{Soit $ A \in \mathbb{P} $ et $ \bar{A} \in \mathbb{P} $ sont complémentaires dans  $ \mathbb{P} $  , soit F une fonction arithmétique complètement multiplicative :}
	\[ \forall F \in \mathbb{M_{c}} ~:~ D(\mathds{1}_{A} \times F,s) \times D(\mathds{1}_{\bar{A}} \times F,s) = D(F,s) \]
\end{proposition}

\begin{preuve}
	\textnormal{Soit $ A \in \mathbb{P} $ et $ \bar{A} \in \mathbb{P} $ sont complémentaire dans  $ \mathbb{P} $  , soit F une fonction arithmétique complètement multiplicative :}
	\begin{align*}
	D(\mathds{1}_{A} \times F,s) \times D(\mathds{1}_{\bar{A}} \times F,s) 
	&= D([\mathds{1}_{A} \times F ] \times [ \mathds{1}_{\bar{A}} \times F ] ,s) \times D([ \mathds{1}_{A} \times F] \scalebox{0.5}{$\square$} [ \mathds{1}_{\bar{A}} \times F],s) \\
	&= D( \mathds{1}_{A} \times F \times  \mathds{1}_{\bar{A}} \times F,s) \times D([ \mathds{1}_{A} \times F] \scalebox{0.5}{$\square$} [ \mathds{1}_{\bar{A}} \times F],s)\\
	&= D( \mathds{1}_{A} \times  \mathds{1}_{\bar{A}} \times F  \times F ,s) \times D([ \mathds{1}_{A} \scalebox{0.5}{$\square$}  \mathds{1}_{\bar{A}}] \times F,s)\\
	&= D( \mathds{1}_{\emptyset} \times F  \times F ,s) \times D(\mathds{1} \times F,s) \\
	&= D( F,s)
	\end{align*}

\end{preuve}

\begin{remark}
	\textnormal{Cette formule reflète simplement le produit d'Euler , prenons zêta comme exemple :}
	\[ D(\mathds{1},s) = \zeta(s) =  \prod_{p \in \mathbb{P}} \frac{1}{1-\frac{1}{p^{s}}} =  \prod_{p \in A} \frac{1}{1-\frac{1}{p^{s}}} \times \prod_{p \in \bar{A}} \frac{1}{1-\frac{1}{p^{s}}} = D(\mathds{1}_{A} \times \mathds{1},s) \times D(\mathds{1}_{\bar{A}} \times \mathds{1} ,s)  \]
\end{remark}

\begin{proposition}
	
	\textnormal{Identité dans $ (\mathbb{M},  \scalebox{0.6}{$\square$}  , \times )  $ : }
	\[[F ~ \scalebox{0.6}{$\square$} ~ F ] (m) = [{2}^{\omega(.)} \times F ](m) \]

\end{proposition}

\begin{preuve}
	\[[F ~ \scalebox{0.6}{$\square$} ~ F ] (m) = \prod_{p|m} F({p}^{v_{p}(m)})+F({p}^{v_{p}(m)}) = [{2}^{\omega(.)} \times F ](m) \]
\end{preuve}

\begin{proposition}
	\textnormal{Identité dans $ (\mathbb{M},  \scalebox{0.6}{$\square$}  , \times )  $ : }
	\[\scalebox{1}{$ \text{Id}_{e} $}  \scalebox{0.6}{$\square$} \mathds{1} (m) = \hat{\sigma}(m)\]
	\textnormal{ici sigma est la somme des diviseurs premiers avec leurs complémentaires .}
\end{proposition}

\begin{preuve}
	
	\[\scalebox{1}{$ \text{Id}_{e} $}  \scalebox{0.6}{$\square$} \mathds{1} (m) = \prod_{p|m} p^{v_{p}(m)} + 1 = \sum_{ab=m ~ a \wedge b = 1} a = \hat{\sigma}(m)\]
\end{preuve}

\begin{proposition}
	\textnormal{Identité dans $ (\mathbb{M},  \scalebox{0.6}{$\square$}  , \times )  $ : }
	\[\phi = \scalebox{1}{$ \text{Id}_{e} $} \times \Bigg[\mathds{1}\scalebox{0.6}{$\square$}\frac{(-1)^{\omega}}{\mathbf{rad}}\Bigg]\]
\end{proposition}

\begin{preuve} 
	La réécriture de la fonction $ \phi $ se déduit simplement de la formule suivante :
	\[\phi(n) = n \times \prod_{i=1}^{\omega(n)}{\left(1-\frac1{p_i}\right)} \]
\end{preuve}

\begin{proposition}

	\[
	\forall s \in \mathbb{C} , \Re(s)>1 ~ : ~
	\Bigg[
	\sum_{p \in \mathbb{P}} \frac{1}{p^{s}} 
	\Bigg] 
	\times 
	\Bigg[ 
	\sum_{n=1}^{\infty} \frac{1}{n^{s}}  
	\Bigg] = \sum_{n=1}^{\infty} \frac{\omega(n)}{n^{s}}
	\]
\end{proposition} 

\begin{preuve}
	
	\begin{align*}
		\Bigg[
		\sum_{p \in \mathbb{P}} \frac{1}{p^{s}} 
		\Bigg] 
		\times 
		\Bigg[ 
		\sum_{n=1}^{\infty} \frac{1}{n^{s}}  
		\Bigg] 
		&= \sum_{p \in \mathbb{P}} \sum_{n=1}^{\infty} \frac{1}{p^{s}} \cdot \frac{1}{n^{s}} \\
		&=\sum_{m=2}^{\infty} \sum_{\underset{p \in \mathbb{P} ~ n \in \mathbb{N^\star}}{pn=m} ~ }\frac{1}{p^{s}} \cdot \frac{1}{n^{s}} \\
		&=\sum_{m=2}^{\infty} \sum_{\underset{p \in \mathbb{P} ~ n \in \mathbb{N^\star}}{pn=m} ~ }\frac{1}{m^{s}} \\
		&=\sum_{m=2}^{\infty} \frac{1}{m^{s}} \sum_{\underset{p \in \mathbb{P} ~ n \in \mathbb{N^\star}}{pn=m} ~ }1 \\
		&= \sum_{m=2}^{\infty} \frac{1}{m^{s}} \cdot \omega(m) \\
		&= \sum_{m=1}^{\infty} \frac{\omega(m)}{m^{s}}
	\end{align*}
	
\end{preuve}

\begin{proposition}
	\[
	\forall s \in \mathbb{C} , \Re(s)>1 ~ : ~
	\sum_{n=1}^{\infty} \frac{1}{n^{s}} \cdot \sum_{p \in \mathbb{P}} \frac{1}{p^{s}} 
	\frac{1}{\omega(np)}  = \zeta(s)-1
	\]

\end{proposition}

\begin{preuve}
	\begin{align*}
		\sum_{n=1}^{\infty} \frac{1}{n^{s}} \cdot \sum_{p \in \mathbb{P}} \frac{1}{p^{s}} 
		 \frac{1}{\omega(np)} 
		&= \sum_{n=1}^{\infty} \sum_{p \in \mathbb{P}}  \frac{1}{p^{s}} \cdot \frac{1}{n^{s}} \cdot \frac{1}{\omega(np)}  \\
		&=\sum_{m=2}^{\infty} \sum_{\underset{p \in \mathbb{P} ~ n \in \mathbb{N^\star}}{pn=m} ~ }\frac{1}{np^{s}} \cdot \frac{1}{\omega(np)} \\
		&=\sum_{m=2}^{\infty} \sum_{\underset{p \in \mathbb{P} ~ n \in \mathbb{N^\star}}{pn=m} ~ }\frac{1}{m^{s}} \cdot \frac{1}{\omega(m)} \\
		&=\sum_{m=2}^{\infty} \frac{1}{m^{s}} \cdot \frac{1}{\omega(m)} \sum_{\underset{p \in \mathbb{P} ~ n \in \mathbb{N^\star}}{pn=m} ~ } 1 \\
		&=\sum_{m=2}^{\infty} \frac{1}{m^{s}} \cdot \frac{1}{\omega(m)} \cdot \omega(m) \\
		&=\sum_{m=2}^{\infty} \frac{1}{m^{s}} \\
		&=\zeta(s)-1 
	\end{align*}
	
\end{preuve}

\begin{proposition}
	\[
	\forall s \in \mathbb{C} , \Re(s)>1 ~ : ~
	\sum_{n=1}^{\infty} \frac{1}{n^{s}} \cdot \frac{1}{\omega(n)+1} \cdot [ \sum_{p \in \mathbb{P}} \frac{1}{p^{s}} +  \sum_{p | n} \frac{1}{p^{s}} \cdot \frac{1}{\omega(n)} ]  = \zeta(s)-1
	\]
\end{proposition}

\begin{preuve}

	\begin{align*}
		\zeta(s)-1
		&= \sum_{n=1}^{\infty} \frac{1}{n^{s}} \cdot \sum_{p \in \mathbb{P}} \frac{1}{p^{s}} 
		\frac{1}{\omega(np)}\\
		&= \sum_{n=1}^{\infty} \sum_{p | n} \frac{1}{p^{s}} \cdot  \frac{1}{\omega(n)} \cdot \frac{1}{n^{s}} + \frac{1}{\omega(n)+1} \cdot \sum_{p \nmid n} \frac{1}{p^{s}} \cdot \frac{1}{n^{s}} \\ 
		&= \sum_{n=1}^{\infty} \frac{1}{n^{s}} \cdot [\sum_{p | n} \frac{1}{p^{s}} \cdot \frac{1}{\omega(n)} + \frac{1}{\omega(n)+1} \cdot \sum_{p \in \mathbb{P}} \frac{1}{p^{s}}- \sum_{p | n} \frac{1}{p^{s}} \cdot  \frac{1}{\omega(n)+1}  ]  \\ 
		&= \sum_{n=1}^{\infty} \frac{1}{n^{s}} \cdot [ \frac{1}{\omega(n)+1} \cdot \sum_{p \in \mathbb{P}} \frac{1}{p^{s}} + \sum_{p | n} \frac{1}{p^{s}}  \cdot [ \frac{1}{\omega(n)}  - \frac{1}{\omega(n)+1} ] ] \\ 
		&= \sum_{n=1}^{\infty} \frac{1}{n^{s}} \cdot [ \frac{1}{\omega(n)+1} \cdot \sum_{p \in \mathbb{P}} \frac{1}{p^{s}} + \sum_{p | n} \frac{1}{p^{s}}  \cdot [ \frac{1}{\omega(n)} \cdot \frac{1}{\omega(n)+1} ] ] \\ 
		&= \sum_{n=1}^{\infty} \frac{1}{n^{s}} \cdot [ \frac{1}{\omega(n)+1} \cdot \sum_{p \in \mathbb{P}} \frac{1}{p^{s}} + \frac{1}{\omega(n)+1} \cdot \sum_{p | n} \frac{1}{p^{s}}  \cdot  \frac{1}{\omega(n)} ] \\ 
		&= \sum_{n=1}^{\infty} \frac{1}{n^{s}} \cdot \frac{1}{\omega(n)+1} \cdot [ \sum_{p \in \mathbb{P}} \frac{1}{p^{s}} + \sum_{p | n} \frac{1}{p^{s}} \cdot \frac{1}{\omega(n)} ]  
	\end{align*}
\end{preuve}

\begin{defini}
	\textnormal{Pour $ y \in \mathbb{R} $ : }
	
	\textnormal{Notons les fonctions $ n \in \mathbb{N^{\star}} ~ : ~  n \rightarrow \cos(y\ln(n)) $  et $ n \in \mathbb{N^{\star}} ~ : ~  n \rightarrow \sin(y\ln(n)) $  respectivement   $\textnormal{Cosa}_{y} \textnormal{ et } \textnormal{Sina}_{y} ~  $ .}
\end{defini}

\begin{proposition}	
	\[\forall m \in \mathbb{N}^{\star}, \forall y \in \mathbb{R} : m^{iy} = [\textnormal{Cosa}_{y} ~ \scalebox{0.6}{$\square$} ~ {i}^{\omega( )}\textnormal{Sina}_{y}](m) \]
\end{proposition}

\begin{preuve}
	\begin{align*}
		\forall n \in \mathbb{N}^{\star} ~ \forall p \in \mathbb{P}, \forall y \in \mathbb{R} : {{p}^{n}}^{iy} 
		&= \exp(iyn\ln(p)) \\
		&= \cos(yn\ln(p))+i\sin(yn\ln(p)) \\
		&= \cos(y\ln(p^{n}))+i\sin(y\ln(p^{n})) \\
		&= [\textnormal{Cosa}_{y} ~ \scalebox{0.6}{$\square$} ~ {i}^{\omega( )}\textnormal{Sina}_{y}]({{p}^{n}}^{iy} ) 
	\end{align*}
\end{preuve}

\begin{proposition}
	\[ 
	\mathds{1} = \textnormal{Cosa}_{y}^{2} ~ \scalebox{0.6}{$\square$} ~ \textnormal{Sina}_{y}^{2} 
	\]

\end{proposition}

\begin{preuve}
	
\begin{align*}
	\forall n \in \mathbb{N}^{\star} ~ \forall p \in \mathbb{P}, \forall y \in \mathbb{R} : \textnormal{Cosa}_{y}^{2} ~ \scalebox{0.6}{$\square$} ~ \textnormal{Sina}_{y}^{2} [{p}^{n}]
	&= \textnormal{Cosa}_{y}^{2}[{p}^{n}] + \textnormal{Sina}_{y}^{2} [{p}^{n}] \\
	&= \cos(y\ln(p^{n}))^{2}+\sin(y\ln(p^{n}))^{2} \\
	&= 1 
	\end{align*}
\end{preuve}

\begin{proposition}
	\textnormal{Identité remarquable : }
	\[A^{2} ~ \scalebox{0.6}{$\square$} ~ (-1)^{\omega(.)}B^{2} = [A ~ \scalebox{0.6}{$\square$} ~ (-1)^{\omega(.)}B ] \times [A ~ \scalebox{0.6}{$\square$} ~ B] \]
\end{proposition}

\begin{preuve}
	\begin{align*}
	\forall A , B \in \mathbb{M} ~ \forall p \in \mathbb{P} ~ \forall n \in \mathbb{N}^{\star}
	\Big[A^{2} ~ \scalebox{0.6}{$\square$} ~ (-1)^{\omega(.)}B^{2} \Big] [p^{n}] 
	&= A(p^{n})^{2}-B(p^{n})^{2} \\
	&= [A(p^{n})-B(p^{n})]\cdot[A(p^{n})+B(p^{n})] \\ 
	&= \Big[[A ~ \scalebox{0.6}{$\square$} ~ (-1)^{\omega(.)}B ] \times [A ~ \scalebox{0.6}{$\square$} ~ B] \Big] [p^{n}]
	\end{align*}
	
\end{preuve}

\begin{proposition}
	\textnormal{les fonctions  $\textnormal{Cosa}_{y}$ et $\textnormal{Sina}_{y}$ vérifient plusieurs identités pour tout $ y \in \mathbb{R^{\star}} $ fixé , des identités qui sont de simple reproduction des formules cos et sin standard dont le seul intérêt est d'être appliqué sur des produits et fractions de nombre  au lieu de sommes et division :}
	
	$$ \textnormal{Cosa}_{y}(ab) = \textnormal{Cosa}_{y}(a)\textnormal{Cosa}_{y}(b)-\textnormal{Sina}_{y}(a)\textnormal{Sina}_{y}(b) $$
	$$ \textnormal{Cosa}_{y}(\frac{a}{b}) = \textnormal{Cosa}_{y}(a)\textnormal{Cosa}_{y}(b)+\textnormal{Sina}_{y}(a)\textnormal{Sina}_{y}(b) $$
	$$ \textnormal{Sina}_{y}(ab) = \textnormal{Sina}_{y}(a)\textnormal{Cosa}_{y}(b)+\textnormal{Cosa}_{y}(a)\textnormal{Sina}_{y}(b) $$
	$$ \textnormal{Sina}_{y}(\frac{a}{b}) = \textnormal{Sina}_{y}(a)\textnormal{Cosa}_{y}(b)-\textnormal{Cosa}_{y}(a)\textnormal{Sina}_{y}(b) $$
	$$ \textnormal{Cosa}_{y}(1) = 1 $$
	$$ \textnormal{Sina}_{y}(1) = 0 $$
	$$ \textnormal{Cosa}_{y}(\frac{1}{a}) = \textnormal{Cosa}_{y}(a) $$
	$$ \textnormal{Sina}_{y}(\frac{1}{a}) = -\textnormal{Sina}_{y}(a) $$
	$$ \textnormal{Cosa}_{2y}(a) + \textnormal{Cosa}_{2y}(b) = 2.\textnormal{Cosa}_{y}(ab).\textnormal{Cosa}_{y}(\frac{a}{b}) $$
	$$ \textnormal{Sina}_{2y}(a) + \textnormal{Sina}_{2y}(b) = 2.\textnormal{Sina}_{y}(ab).\textnormal{Cosa}_{y}(\frac{a}{b}) $$
	$$ \textnormal{Sina}_{2y}(a) = 2.\textnormal{Cosa}_{y}(a)\textnormal{Sina}_{y}(a) $$
	$$ \textnormal{Cosa}_{y}(a)^{2} + \textnormal{Cosa}_{y}(b)^{2} -\textnormal{Cosa}_{y}(ab).\textnormal{Cosa}_{y}(\frac{a}{b}) = 1 $$
	$$ \textnormal{Cosa}_{y}(a^{2}) = \textnormal{Cosa}_{2y}(a)   $$
	$$ \textnormal{Sina}_{y}(a^{2}) = \textnormal{Sina}_{2y}(a)   $$
	$$ \textnormal{Cosa}_{y}(a).\textnormal{Cosa}_{y}(b) =\frac{\textnormal{Cosa}_{y}(ab)+\textnormal{Cosa}_{y}(\frac{a}{b})}{2} $$
	$$ \textnormal{Sina}_{y}(a).\textnormal{Sina}_{y}(b) =\frac{\textnormal{Cosa}_{y}(\frac{a}{b})-\textnormal{Cosa}_{y}(ab)}{2} $$
	$$ \textnormal{Cosa}_{y}(a).\textnormal{Cosa}_{y}(b).\textnormal{Cosa}_{y}(c) =\frac{\textnormal{Cosa}_{y}(abc)+\textnormal{Cosa}_{y}(\frac{ab}{c})+\textnormal{Cosa}_{y}(\frac{ac}{b})+\textnormal{Cosa}_{y}(\frac{bc}{a})}{4} $$
	
\end{proposition}

\begin{remark}
	\textnormal{la fonction utilisée dans le premier chapitre peut maintenant s'écrire sous une autre forme :}
	\begin{align*}
		 Q(l) 
		 &= \Bigg[ \delta_{l=1}(l) + \sum_{\underset{mk=l, m>k }{m \wedge k = 1} } 2\Re\Big(\frac{\chi(m)}{{m}^{iy}}\overline{\frac{\chi(k)}{{k}^{iy}}}\Big) \Bigg]\\
		 &=[\textnormal{Cosa}_{y}\times\Re\chi ~ \scalebox{0.6}{$\square$} ~ \textnormal{Sina}_{y}\times\Im\chi] [l]
	\end{align*}
	
	Avec :
	\[ \forall n \in \mathbb{N}^{\star} ~ \forall p \in \mathbb{P} ~: ~ \\
	\Re\chi [p^{n}] =  \Re(\chi(p^{v_{p}(n)}))\]
	\[ \forall n \in \mathbb{N}^{\star} ~ \forall p \in \mathbb{P} ~: ~ \\
	\Im\chi [p^{n}] =  \Im(\chi(p^{v_{p}(n)}))\]
\end{remark}

\begin{remark}
	On définit les trois fonctions suivantes : 
	
	\[ \forall n \in \mathbb{N}-\{0,1\} , ~ \forall a \in \mathbb{N}-\{0,1\} \textrm{ : } \bigg[ \stackrel{\phi(n)}{  \underset{k=1}{  \mathlarger{ \mathlarger{ \mathlarger{ \square } } } } }    \mathlarger{ \mathlarger{ \chi_{k} } }   \bigg] (a) \]
	
	\[ \forall l \in \mathbb{Z}/n\mathbb{Z} : \prod_{k=1}^{\phi(n)} \mathlarger{ \mathlarger{ \chi_{k} } } (l) \]
	
	\[ \forall n \in \mathbb{N}-\{0,1\} : \sum_{l \in \mathbb{Z}/n\mathbb{Z}}{\prod_{k=1}^{\phi(n)} \mathlarger{ \mathlarger{ \chi_{k} }} (l)} \]

\end{remark}

\begin{remark}
	Les fonctions présentées ci dessous sont simplement citées car il pourraient être utilisé (notamment la première formule ) pour calculer des suites dépendantes de tous les caractères du même entier.
	
	Les deux dernières fonctions présentes numériquement les observations suivantes :

	\[ \forall l \in \mathbb{Z}/n\mathbb{Z} : \prod_{k=1}^{\phi(n)} \raisebox{2pt}{ $\chi_{k}$} (l) = \begin{cases}
	\begin{array}{ll}
	0 \, \textnormal{ si } l \wedge n > 1\\
	\pm1 \, \textnormal{ si non} \\
	\end{array}
	\end{cases} =  \begin{cases}
	\begin{array}{ll}
	\chi_{\frac{\phi(n)}{2}}(l) \,\textnormal{ si ? }\\ -\chi_{\frac{\phi(n)}{2}}(l) \,\textnormal{ si ? } 
	\end{array}
	\end{cases}	\]
	
	\[ \forall n \in \mathbb{N}-\{0,1\} : \sum_{l \in \mathbb{Z}/n\mathbb{Z}}{\prod_{k=1}^{\phi(n)} \raisebox{2pt}{ $ \chi_{k} $ }} (l)=  \begin{cases}
	\begin{array}{ll}
	0 \,\textnormal{ si } \,\textnormal{  ?} \\
	\phi(n) \,\textnormal{ ?}
	\end{array}
	\end{cases} 	\]

\end{remark}

\begin{proposition} 
	\[ \forall n \in \mathbb{N}-\{0,1\} ,~ \forall a \in \mathbb{N}-\{0,1\}  :  \bigg[ \stackrel{\phi(n)} { \mathlarger{ \mathlarger{ \underset{k=1}{ ~  \mathlarger{ \mathlarger{ \square } }  ~ } }}  }   \mathlarger{ \mathlarger{ \chi_{k} } }  \bigg] (a)  = \phi(n)^{\omega(a)} \cdot \delta_{n | [\text{Id}_{e} - \mathds{1}](a) } \]
\end{proposition}

\begin{preuve}
\begin{align*}
	 \forall n \in \mathbb{N}-\{0,1\} ~ \forall a \in \mathbb{N}-\{0,1\}  :  \bigg[\stackrel{\phi(n)} {  \mathlarger{ \underset{k=1}{  \mathlarger{ \mathlarger{ \square } }  } }   }\chi_{k}  \bigg] (a)
	 &=\prod_{p|a} \sum_{k=1}^{\phi(n)} \chi_{k}({p}^{v_{p}(a)}) \\
	 &=\prod_{p|a} \phi(n) \cdot \delta_{{p}^{v_{p}(a)} \equiv 1 \mod n} \\
	 &=\prod_{p|a} \phi(n) \cdot \delta_{{p}^{v_{p}(a)}-1 \equiv 0 \mod n} \\
	 &=\phi(n)^{\omega(a)} \cdot \delta_{[\prod_{p|a} {p}^{v_{p}(a)}-1] \equiv 0 \mod n} \\
	 &=\phi(n)^{\omega(a)} \cdot \delta_{[ \text{Id}_{e} -\mathds{1}](a) \equiv 0 \mod n} \\
	 &=\phi(n)^{\omega(a)} \cdot \delta_{n | [  \text{Id}_{e} - \mathds{1}](a) } \\
\end{align*}
\end{preuve}

\iffalse

\begin{remark}
	%\mathlarger{ \mathlarger{ \chi_{k} } }
	Je n’ai pas réussi à définir une dérivée dans l'anneau construit précédemment , une des idées derrière la construction de cette convolution est de définir une dérivée de fonction multiplicative ( ou peut être simplement des fonctions complètement multiplicatives ) .\\
	
	l'idée est d'amortir les oscillations des séries qui convergent au sens de Cesàro sans converger. la solution parfaite est de trouver pour chaque fonction multiplicative dont la série de Dirichlet ne converge pas mais qui converge au sens de Cesàro , une fonction multiplicative qui convergerait vers la limite de Cesàro .\\
	
	%Un exemple explicatif serait de voir si des fonctions du style  $ F_{amortie} := \int \frac{F'}{n} $ ont des séries de Dirichlet convergentes vers la limite de Cesàro de la série $ D(F,s) $ , ceci est loin d'être ne serait ce que possible , je ne donne cette exemple que pour illustrer la démarche .	
	Un problème intéressant serait d'essayer de donner un sens adapté au dérivée de fonctions multiplicatives.
	
	Une telle dérivée devrait vérifier les contraintes suivantes :
	
	\[ (F {\, \scalebox{0.5}{$\square$} \,} G)' = F' {\, \scalebox{0.5}{$\square$} \,} G'\]
	\[ (F \times G)' = F' \times G {\, \scalebox{0.5}{$\square$} \,} F \times G' \]
	une approche serait de se restreindre au fonction multiplicative polynomial en les nombres premiers , ou en les puissances de nombres premiers . Ou alors définir la dérivée uniquement pour des fonctions complètement multiplicatives.
\end{remark} 
\fi
\begin{remark}
	il serait intéressant de pouvoir amortir les oscillations des séries qui convergent au sens de Cesàro sans converger. la solution parfaite est de trouver pour chaque fonction multiplicative dont la série de Dirichlet ne converge pas mais qui converge au sens de Cesàro , une fonction multiplicative qui convergerait vers la limite de Cesàro .\\
\end{remark}